\documentclass[11pt]{article}
\usepackage{amssymb,amsmath}
\usepackage{arrows}
\numberwithin{equation}{section}
\allowdisplaybreaks[1]
\newcommand{\real}{{\mathbb R}}
\newcommand{\complex}{{\mathbb C}}
\renewcommand{\natural}{{\mathbb N}}

\renewcommand{\Re}{\mathrm{Re}}

\newcommand{\Ker}{\mathrm{Ker}}
\newcommand{\Ran}{\mathrm{Ran}}
\renewcommand{\div}{\mathop{\rm div}}

\renewcommand{\d}{\,{\rm d}}            
\newcommand{\D}{{\rm d}}                

\newcommand{\weakto}{\rightharpoonup}

\newcommand{\cE}{{\cal E}}
\newcommand{\cF}{{\cal F}}

\newcommand{\cI}{{\cal I}}
\newcommand{\cJ}{{\cal J}}
\newcommand{\cK}{{\cal K}}
\newcommand{\cL}{{\cal L}}
\newcommand{\cM}{{\cal M}}

\newcommand{\cO}{{\cal O}}

\newcommand{\cS}{{\cal S}}

\newtheorem{theorem}{Theorem}[section]
\newtheorem{lemma}[theorem]{Lemma}

\newtheorem{proposition}[theorem]{Proposition}

\newtheorem{remark}[theorem]{Remark}

\newcommand{\proof}{{\noindent \bf Proof.\ }}

\newcommand{\DS}{\displaystyle}
\newcommand{\app}{{\rm app}}

\newcommand{\tv}{{\rm tv}}
\newcommand{\inttwo}{{\int_{\real^2}}}

\def\build#1_#2^#3{\mathrel{
  \mathop{\kern 0pt#1}\limits_{#2}^{#3}}}

\def\QED{\mbox{}\hfill$\Box$}

\newdimen\texpscorrection
\newdimen\figcenter
\def\figurewithtex #1 #2 #3 #4 #5\cr{\null
  {\goodbreak\figcenter=\hsize\relax
  \advance\figcenter by -#4truecm
  \divide\figcenter by 2
  \begin{figure}[hbt]
  \vskip #3truecm\noindent\hskip\figcenter
  \includegraphics{#1}{\hskip\texpscorrection\input #2 }
  \vskip 0.8truecm{\baselineskip=0.8\baselineskip
  \noindent \vbox{\noindent {\footnotesize #5}}\par}
  \end{figure}}}
\def\point#1 #2 #3 {\rlap{\kern #1 truecm
\raise #2 truecm \hbox{#3}}}

\begin{document}

\title{Interaction of vortices in viscous planar flows}

\author{{\bf Thierry Gallay} \\[1mm] 
Universit\'e de Grenoble I\\
Institut Fourier, UMR CNRS 5582\\
B.P. 74\\
F-38402 Saint-Martin-d'H\`eres, France\\
{\tt Thierry.Gallay@ujf-grenoble.fr}}
\date{August 18, 2009}

\maketitle

\begin{abstract}
We consider the inviscid limit for the two-dimensional
incompressible Navier-Stokes equation in the particular case where
the initial flow is a finite collection of point vortices. We
suppose that the initial positions and the circulations of the
vortices do not depend on the viscosity parameter $\nu$, and we
choose a time $T > 0$ such that the Helmholtz-Kirchhoff point vortex
system is well-posed on the interval $[0,T]$. Under these
assumptions, we prove that the solution of the Navier-Stokes
equation converges, as $\nu \to 0$, to a superposition of Lamb-Oseen
vortices whose centers evolve according to a viscous regularization
of the point vortex system. Convergence holds uniformly in time, in
a strong topology which allows to give an accurate description of
the asymptotic profile of each individual vortex. In particular, we
compute to leading order the deformations of the vortices due to
mutual interactions. This allows to estimate the self-interactions, 
which play an important role in the convergence proof.
\end{abstract}

\section{Introduction} \label{sec1}

It is a well established fact that coherent structures play a crucial
role in the dynamics of two-dimensional turbulent flows. Experimental
observations \cite{Cou83} and numerical simulations of decaying
turbulence \cite{MW84,MW90} reveal that, in a two-dimensional flow
with sufficiently high Reynolds number, isolated regions of
concentrated vorticity appear after a short transient period, and
persist over a very long time scale. These structures are nearly
axisymmetric and behave like point vortices as long as they remain
widely separated, but when two of them come sufficiently close to each
other they get significantly deformed under the the strain of the
velocity field, and the interaction may even cause both vortices to
merge into a single, larger structure \cite{LD02,Meu05}. It thus
appears that the long-time behavior of two-dimensional decaying
turbulence is essentially governed by a few basic mechanisms, such as
vortex interaction and, especially, vortex merging.

Although these phenomena are relatively well understood from a
qualitative point of view, they remain largely beyond the scope of
rigorous analysis. Vortex merging, in particular, is a genuinely
nonperturbative process which seems extremely hard to describe
mathematically, although it is certainly the key mechanism which
explains the coarsening of vorticity structures in two-dimensional
flows, in agreement with the inverse energy cascade. The situation is
better for vortex interactions, which may be rigorously studied in the
asymptotic regime where the distance between vortices is much larger
than the typical core size, but complex phenomena can occur even in
that case. Indeed, numerical calculations \cite{LD02} and nonrigorous
asymptotic expansions \cite{TT65,TK91} indicate that vortex
interaction begins with a fast relaxation process, during which each
vortex adapts its shape to the velocity field generated by the other
vortices. This first step depends on the details of the initial data,
and is characterized by temporal oscillations of the vortex cores
which disappear on a non-viscous time scale. In a second step, the
vortices relax to a Gaussian-like profile at a diffusive rate, and the
system reaches a ``metastable state'' which is independent of the
initial data, and will persist until two vortices get sufficiently
close to start a merging process. In this metastable regime, the
vortex centers move in the plane according to the Helmholtz-Kirchhoff
dynamics, and the vortex profiles are uniquely determined, up to a
scaling factor, by the relative positions of the centers.

From a mathematical point of view, a natural approach to study vortex
interactions is to start with {\em point vortices} as initial data.
After solving the Navier-Stokes equations, we obtain in this way a family
of interacting vortices which, by construction, is directly in the
metastable state that we have just described. In particular, as it
will be proved below, we do not observe here the oscillatory and
diffusive transient steps which take place in the general case. Point
vortices can therefore be considered as {\em well prepared initial
data} for the vortex interaction problem.

With this motivation in mind, we study in the present paper what we
call the {\em viscous $N$-vortex solution}, namely the solution of the
two-dimensional Navier-Stokes equations in the particular case where
the initial vorticity is a superposition of $N$ point vortices. For a
given value of the viscosity parameter $\nu$, this solution is
entirely determined by the initial positions $x_1,\dots,x_N$ and the
circulations $\alpha_1,\dots, \alpha_N$ of the vortices. It describes
a family of interacting vortices of size $\cO((\nu t)^{1/2})$, which
are therefore widely separated if $\nu$ is sufficiently small. Our
main goal is to obtain a rigorous asymptotic expansion of the
$N$-vortex solution in the vanishing viscosity limit, assuming that
vortex collisions do not occur. As was already explained, this problem
is physically relevant, but it also has its own mathematical
interest. Indeed, it is known that the two-dimensional Navier-Stokes
equations have a unique solution, for any value of $\nu$, when the
initial vorticity is a finite measure \cite{GG05}, but computing the
inviscid limit of rough solutions is a very difficult task in general,
due to the underlying instabilites of the Euler flow. Surprisingly
enough, although point vortices are perhaps the most singular initial
data that can be allowed for the Navier-Stokes equations, the inviscid
limit appears to be tractable for the $N$-vortex solution, and
provides a new rigorous derivation of the Helmholtz-Kirchhoff dynamics
as well as a mathematical description of the metastable regime for
interacting vortices.

In the rest of this introductory section, we recall a global
well-posedness result for the two-dimensional Navier-Stokes equations
which is adapted to our purposes, we introduce the Lamb-Oseen vortices
which will play a crucial role in our analysis, and we briefly mention
the difficulties related to the inviscid limit of rough solutions. Our
main results concerning the $N$-vortex solution will be stated in
Section~\ref{sec2}, and proved in the subsequent sections.

The incompressible Navier-Stokes equations in the plane $\real^2$ 
have the following form:
\begin{equation}\label{NS}
  \frac{\partial u}{\partial t} + (u \cdot \nabla)u \,=\, 
  \nu \Delta u - \nabla p~, \qquad \div u = 0~,
\end{equation}
where $u(x,t) \in \real^2$ denotes the velocity of the fluid at point $x
\in \real^2$ and time $t > 0$, and $p(x,t) \in \real$ is the pressure
inside the fluid. The only physical parameter in \eqref{NS} is the
{\em kinematic viscosity} $\nu > 0$, which will play an important role
in this work. For our purposes it will be convenient to consider the
{\em vorticity field} $\omega(x,t) =
\partial_1 u_2(x,t) - \partial_2 u_1(x,t)$, which evolves according 
to the remarkably simple equation
\begin{equation}\label{V}
  \frac{\partial\omega}{\partial t} + (u \cdot \nabla) \omega \,=\, 
  \nu \Delta \omega~. 
\end{equation}
Under mild assumptions, which will always be satisfied below, 
the velocity field $u(x,t)$ can be reconstructed from the vorticity 
$\omega(x,t)$ via the two-dimensional Biot-Savart law:
\begin{equation}\label{BS}
  u(x,t) \,=\, \frac{1}{2\pi} \inttwo 
  \frac{(x - y)^{\perp}}{|x -y|^2} \,\omega(y,t)\d y\ ,
  \quad x \in \real^2~,
\end{equation}
where, for any $x = (x_1,x_2) \in \real^2$, we denote $x^\perp = (-x_2,x_1)$
and $|x|^2 = x_1^2 + x_2^2 $. 

Let $\cM(\real^2)$ be the space of all real-valued finite measures 
on $\real^2$, equipped with the total variation norm
$$
  \|\mu\|_\tv \,=\, \sup\{\langle\mu_n,\phi\rangle \,;\, \phi \in 
  C_0(\real^2)\,,~ \|\phi\|_{L^\infty} \le 1\}~.
$$
Here $\langle\mu,\phi\rangle = \inttwo \phi\d\mu$, and $C_0(\real^2)$
denotes the space of all continuous functions $\phi : \real^2 \to
\real$ which vanish at infinity. We say that a sequence $\{\mu_n\}$ in
$\cM(\real^2)$ converges weakly to $\mu \in \cM(\real^2)$ if
$\langle\mu_n,\phi \rangle \to \langle \mu,\phi\rangle$ as $n \to
\infty$ for all $\phi \in C_0(\real^2)$. Weak convergence is denoted
by $\mu_n \weakto \mu$.

Our starting point is the following result, which shows that 
the initial value problem for Eq.~\eqref{V} is globally 
well-posed in the space $\cM(\real^2)$:

\begin{theorem}\label{thm1} {\bf \cite{GG05}}
Fix $\nu > 0$. For any initial measure $\mu \in \cM(\real^2)$, 
Eq.~\eqref{V} has a unique global solution 
\begin{equation}\label{omegL1}
  \omega \in C^0((0,\infty),L^1(\real^2) \cap 
  L^\infty(\real^2))
\end{equation}
such that $\|\omega(\cdot,t)\|_{L^1} \le \|\mu\|_\tv$ for all 
$t > 0$, and $\omega(\cdot,t) \weakto \mu$ as $t \to 0+$. 
\end{theorem}

Here and in what follows, it is understood that $\omega$ is 
a {\em mild} solution of \eqref{V}, i.e. a solution of the
associated integral equation
\begin{equation}\label{integ}
  \omega(t) \,=\, e^{\nu t\Delta}\mu - \int_0^t \div\Bigl(
  e^{\nu (t-s)\Delta} u(s)\omega(s)\Bigr)\d s~, \quad t > 0~,
\end{equation}
where $e^{t\Delta}$ denotes the heat semigroup. The {\em existence} of
a global solution to \eqref{V} for all initial data in $\cM(\real^2)$
has been established more than 20 years ago by G.-H.~Cottet
\cite{Co96}, and independently by Y.~Giga, T.~Miyakawa and H.~Osada
\cite{GMO88}. In the same spirit, the later work by T.~Kato
\cite{Ka94} should also be mentioned. In addition to existence, it was
shown in \cite{GMO88,Ka94} that the solution of \eqref{V} is {\em
  unique} if the atomic part of the initial measure is small compared
to the viscosity. This smallness condition turns out to be necessary
if one wants to obtain uniqueness by a standard application of
Gronwall's lemma. On the other hand, in the particular case where the
initial vorticity is a single Dirac mass (of arbitrary size),
uniqueness of the solution of \eqref{V} was proved recently by
C.E.~Wayne and the author \cite{GW05}, using a dynamical system
approach. An alternative proof of the same result can also be found in
\cite{GGL05}. Finally, in the general case, it is possible to obtain 
uniqueness of the solution of \eqref{V} by isolating the large Dirac 
masses in the initial measure and combining the approaches of 
\cite{GMO88} and \cite{GW05}. This last step in the proof of 
Theorem~\ref{thm1} was achieved by I.~Gallagher and the author 
in \cite{GG05}.

When the initial vorticity $\mu = \alpha \delta$ is a multiple
of the Dirac mass (located at the origin), the unique solution 
of \eqref{V} is an explicit self-similar solution called the
{\em Lamb-Oseen vortex:} 
\begin{equation}\label{Oseen}
  \omega(x,t) \,=\, \frac{\alpha}{\nu t}\,G\Bigl(\frac{x}
  {\sqrt{\nu t}}\Bigr)~, \qquad u(x,t) \,=\,\frac{\alpha}{\sqrt{\nu t}}
  \,v^G\Bigl(\frac{x}{\sqrt{\nu t}}\Bigr)~,
\end{equation}
where
\begin{equation}\label{Gdef}
  G(\xi) \,=\, \frac{1}{4\pi} \,e^{-|\xi|^2/4}~,\qquad
  v^G(\xi) \,=\, \frac{1}{2\pi}\frac{\xi^\perp}{|\xi|^2}
  \Bigl(1 -  e^{-|\xi|^2/4}\Bigr)~, \qquad \xi \in \real^2~.
\end{equation}
The circulation parameter $\alpha \in \real$ measures the intensity of 
the vortex. It coincides, in this particular case, with the integral 
of the vorticity $\omega$ over the whole plane $\real^2$, a quantity 
which is preserved under evolution. The dimensionless quantity 
$|\alpha|/\nu$ is usually called the {\em circulation Reynolds number}.

According to Theorem~\ref{thm1}, the vorticity equation \eqref{V} has
a unique global solution for any initial measure and any value of the
viscosity parameter. It is very natural to investigate the behavior of
that solution in the vanishing viscosity limit, especially if the
initial data contain non-smooth structures such as point vortices,
vortex sheets, or vortex patches. Indeed, if the viscosity is small,
these structures will persist over a sufficiently long time scale to
be observed and to influence the dynamics of the system. This
question, however, is very difficult in its full generality, because
Eq.~\eqref{V} reduces formally, as $\nu \to 0$, to the Eulerian
vorticity equation
\begin{equation}\label{E}
  \frac{\partial\omega}{\partial t} + (u \cdot \nabla) \omega \,=\, 0~,
\end{equation}
which is not known to be well-posed in such a large space as
$\cM(\real^2)$. If the initial vorticity $\mu \in \cM(\real^2)$
belongs to $H^{-1}(\real^2)$, and if the singular part of $\mu$ has a
definite sign, then Eq.~\eqref{E} has at least a global weak solution
\cite{De91,Ma93}, but this result does not cover the case of point
vortices due to the assumption $\mu \in H^{-1}(\real^2)$. Furthermore,
if we want to prove that the solution of \eqref{E} is unique, we have to
assume that the initial vorticity is bounded \cite{Yu63} or almost bounded
\cite{Yu95,Vi99}.

\medskip
From a more general point of view, it is relatively easy to show that
solutions of the Navier-Stokes equations converge, in the vanishing
viscosity limit, to solutions of the Euler equations if we restrict
ourselves to {\em smooth} solutions in a domain {\em without boundary}
\cite{EM70,Sw71,Ka72,BM81}. The situation is completely different in
the presence of boundaries, especially if one chooses the classical
non-slip boundary conditions for the solutions of the Navier-Stokes
equations. In that case Prandtl boundary layers may form, and the
inviscid limit becomes an extremely difficult question which has been
rigorously treated so far only for well-prepared analytic data
\cite{SC98,Gr00}. But even in the absence of boundaries, hard problems
can occur in the inviscid limit if one considers non-smooth
solutions. Classical examples are listed below, in increasing order of
singularity:

\medskip\noindent{\bf 1)} {\em Vortex patches}. The simplest example
in this category is the case where the vorticity is the characteristic
function of a smooth bounded domain in $\real^2$. The corresponding
velocity field is Lipschitz continuous, or almost Lipschitz if the
boundary of the patch has singularities. The first convergence
results, due to P.~Constantin and J.~Wu \cite{CW95,CW96}, and to
J.-Y.~Chemin \cite{Ch96}, hold for a general class of solutions
including two-dimensional vortex patches. Several improvements have
been subsequently obtained, especially by R.~Danchin \cite{Da97,Da99}, 
H.~Abidi and R.~Danchin \cite{AD04}, T. Hmidi \cite{Hm05,Hm06}, and 
N. Masmoudi \cite{Ma07}. Closer to the spirit of the present work, we 
also quote a recent paper by F.~Sueur \cite{Su08}, where viscous transition
profiles at the boundary of a vortex patch are systematically
constructed, and provide a complete asymptotic expansion of the
solution in powers of $(\nu t)^{1/2}$.

\medskip\noindent{\bf 2)} {\em Vortex sheets}. In a two-dimensional
setting, this is the case where the initial vorticity is concentrated 
on a piece of curve in $\real^2$. The tangential component of the
velocity field is discontinuous along the curve, thus creating a 
shear flow which is responsible for the celebrated Kelvin-Helmholtz 
instability. Due to this underlying instability, the inviscid limit
for vortex sheets is at least as difficult as for Prandtl boundary layers, 
and has been rigorously treated so far only in the case of well-prepared 
analytic data \cite{CW91,CS06}.

\medskip\noindent{\bf 3)} {\em Point vortices}. This is the most
singular example in our list, since here the velocity field is not
even bounded near the vortex centers. However, in the case of point
vortices, the evolution of the inviscid solution is given, at least
formally, by the Helmholtz-Kirchhoff system which does not exhibit any
dynamical instability. Therefore, one can reasonably hope to control
the vanishing viscosity limit in this particular situation. The first
rigorous results in this direction were obtained by C.~Marchioro
\cite{Ma90,Ma98}, and the aim of the present paper is to show how
these results can be complemented to obtain an accurate description of
the slightly viscous solution of \eqref{V} when the initial condition
is a finite collection of point vortices.

After this general introduction, we now give a precise definition 
of the problem we want to study, and we state our main results. 

\section{The viscous $N$-vortex solution}\label{sec2}

Let $N$ be a positive integer. We take $x_1,\dots,\,x_N \in \real^2$ 
such that $x_i \neq x_j$ for $i \neq j$, and we also fix $\alpha_1,
\dots,\,\alpha_N \in \real\setminus\{0\}$. Given any $\nu > 0$, 
we denote by $\omega^\nu(x,t)$, $u^\nu(x,t)$ the unique solution 
of the vorticity equation \eqref{V} with initial data
\begin{equation}\label{mudef}
  \mu \,=\, \sum_{i=1}^N \alpha_i \,\delta(\cdot - x_i)~.
\end{equation}
This initial measure, which does not depend on the viscosity $\nu$, 
describes a superposition of $N$ point vortices of circulations 
$\alpha_1,\dots,\,\alpha_N$ located at the points $x_1,\dots,\,x_N$
in $\real^2$. Our goal is to describe the behavior of the 
vorticity $\omega^\nu(x,t)$ in the vanishing viscosity limit.

The measure $\mu$ is very singular, and we do not know how to 
construct even a weak solution of Euler's equation \eqref{E} with
such initial data. However, the inviscid motion of point vortices 
in the plane has been investigated by many authors, starting 
with H. von Helmholtz \cite{He58} and G.~R. Kirchhoff \cite{Ki76} 
who derived a system of ordinary differential equations describing 
the motion of the vortex centers. It is therefore reasonable to 
expect, in our case, that
\begin{equation}\label{weakapp}
  \omega^\nu(\cdot,t) ~\approx~ \sum_{i=1}^N \alpha_i 
  \,\delta(\cdot - z_i(t))~, \quad \hbox{as }\nu \to 0~,
\end{equation}
where $z(t) = (z_1(t),\dots,\,z_N(t))$ denotes the solution of the 
Helmholtz-Kirchhoff system
\begin{equation}\label{PW}
  z_i'(t) \,=\, \frac{1}{2\pi} \sum_{j\neq i} \alpha_j 
  \,\frac{(z_i(t) - z_j(t))^\perp}{|z_i(t)-z_j(t)|^2}~, 
  \qquad z_i(0) \,=\, x_i~.
\end{equation}

As was already mentioned, the expression in the right-hand side of
\eqref{weakapp} is {\em not} a weak solution of Euler's equation,
because in that case the self-interaction terms in the nonlinearity
$u\cdot\nabla \omega$ are too singular to make sense even as
distributions. However, as was shown by C.~Marchioro and M.~Pulvirenti
\cite{Ma88,MP93}, it is possible to derive system \eqref{PW} from
Euler's equation by a rigorous procedure, which consists in
approximating the Dirac masses in the initial data by small vortex
patches, of size $\epsilon > 0$ and circulations
$\alpha_1,\dots,\,\alpha_N$, whose centers are located at the points
$x_1,\dots,\,x_N$. Then the corresponding solution of \eqref{E}
converges weakly, as $\epsilon \to 0$, to the expression
\eqref{weakapp} where the vortex positions $z_1(t),\dots,\,z_N(t)$ are
solutions of \eqref{PW}. We also recall that system \eqref{PW} is not
globally well-posed for all initial data, because if $N \ge 3$ and if
the circulations $\alpha_i$ are not all the same sign, vortex
collisions can occur in finite time for some exceptional initial
configurations \cite{MP94}.

Our first result shows that the solution $\omega^\nu(x,t)$ 
of Eq.~\eqref{V} given by Theorem~\ref{thm1} converges weakly, 
in the vanishing viscosity limit, to a superposition of 
point vortices which evolve according to \eqref{PW}, provided 
that vortex collisions do not occur.  

\begin{theorem}\label{thm2}
Assume that the point vortex system \eqref{PW} is well-posed 
on the time interval $[0,T]$. Then the solution $\omega^\nu(x,t)$
of the Navier-Sokes equation \eqref{V} with initial data 
\eqref{mudef} satisfies
\begin{equation}\label{weakconv}
  \omega^\nu(\cdot,t) 
  ~\xrightharpoonup[\nu \to 0]{\hbox to 8mm{}}~ 
  \sum_{i=1}^N \alpha_i \,\delta(\cdot - z_i(t))~, \quad 
  \hbox{for all } t \in [0,T]~,
\end{equation}
where $z(t) = (z_1(t),\dots,\,z_N(t))$ is the solution of \eqref{PW}. 
\end{theorem}

This theorem is closely related to a result by C. Marchioro 
\cite{Ma90,Ma98}, which we now briefly describe. Instead of 
point vortices, Marchioro considers initial data of the form
\[
  \omega_0^\epsilon(x) \,=\, \sum_{i=1}^N \omega_i^\epsilon(x)~, 
  \qquad \epsilon > 0~,
\]
where, for each $i \in \{1,\dots,N\}$, $\omega_i^\epsilon$ is a smooth 
vortex patch with a definite sign, which is centered at point 
$x_i \in \real^2$, has compact support of size $\cO(\epsilon)$, 
and satisfies
\[
  \inttwo \omega_i^\epsilon(x)\d x \,=\, \alpha_i^\epsilon 
  ~\xrightarrow[\epsilon \to 0]{}~ \alpha_i~.
\]
Under these assumptions, it is proved that the solution
$\omega^{\epsilon,\nu}(x,t)$ of \eqref{V} with initial data
$\omega_0^\epsilon$ converges to the expression in the right-hand 
side of \eqref{weakapp} in the double limit $\nu \to 0$, 
$\epsilon \to 0$, provided
\begin{equation}\label{epsnu}
  \nu \,\le\, \nu_0\,\epsilon^\beta~, \quad \hbox{for some } \nu_0 > 0
  \hbox{ and } \beta > 0~.
\end{equation}
Theorem~\ref{thm2} above corresponds to the limiting case
$\epsilon = 0$, $\nu \to 0$, which is precisely excluded by 
hypothesis~\eqref{epsnu}. It should be mentioned, however, 
that restriction \eqref{epsnu} can be removed if the circulations 
$\alpha_1,\dots,\alpha_N$ are all the same sign, in which case 
Theorem~\ref{thm2} may probably be established using the techniques
developped in \cite{Ma90}. 

Marchioro's proof is based on a decomposition of the solution
$\omega^{\epsilon,\nu}(x,t)$ into a sum of $N$ viscous vortex
patches. The main idea is to control the spread of each patch by
computing its moment of inertia with respect to a suitable point
$z_i^{\epsilon,\nu}(t)$, which is an approximate solution of
\eqref{PW}. This argument does not give any information on the actual
shape of the vortex patches, and is therefore not sufficient to
provide a qualitative description of the solution
$\omega^{\epsilon,\nu}(x,t)$ for small $\epsilon,\nu$.
Theorem~\ref{thm2} above suffers exactly from the same drawback.  Its
main interest is to provide a natural and rigorous derivation of the
point vortex system \eqref{PW}, which differs from the classical
approach of \cite{Ma88,MP93}.

\medskip The main goal of the present paper is to obtain a
quantitative version of Theorem~\ref{thm2} which specifies the
convergence rate in \eqref{weakconv} and provides a precise asymptotic
expansion of the $N$-vortex solution $\omega^\nu(x,t)$ in the
vanishing viscosity limit. As in Marchioro's approach, our starting
point is a decomposition of $\omega^\nu(x,t)$ into a sum of $N$
viscous vortex patches, which correspond to the atoms of the initial
measure \eqref{mudef}.

\begin{lemma}\label{lemdec}
The $N$-vortex solution can be decomposed as 
\begin{equation}\label{omdecomp}
  \omega^\nu(x,t) \,=\, \sum_{i=1}^N \omega_i^\nu(x,t)~, \qquad
  u^\nu(x,t) \,=\, \sum_{i=1}^N u_i^\nu(x,t)~,  
\end{equation}
where, for each $i \in \{1,\dots,N\}$, $\omega_i^\nu \in C^0((0,\infty),
L^1(\real^2)\cap L^\infty(\real^2))$ is the (unique) solution of the 
convection-diffusion equation
\begin{equation}\label{convecdiff}
  \frac{\partial\omega_i^\nu}{\partial t} + (u^\nu \cdot \nabla) 
  \omega_i^\nu \,=\, \nu \Delta \omega_i^\nu~, \quad \hbox{with}\quad 
  \omega_i^\nu(\cdot,t) ~\xrightharpoonup[t \to 0]{\hbox to 6mm{}}~ 
  \alpha_i \,\delta(\cdot - x_i)~, 
\end{equation}
and the velocity field $u_i^\nu$ is obtained from $\omega_i^\nu$ via 
the Biot-Savart law \eqref{BS}. Moreover $\omega_i^\nu(x,t)$ has the 
same sign as $\alpha_i$ for all $x \in \real^2$ and all $t > 0$, 
and $\inttwo \omega_i^\nu(x,t)\d x = \alpha_i$ for all $t > 0$. 
Finally, there exists $K_0 > 0$ (depending only on $\nu$ and 
$|\alpha| = |\alpha_1| + \dots + |\alpha_N|$) such that
\begin{equation}\label{omegaibdd}
   |\omega_i^\nu(x,t)| \,\le\, K_0\,\frac{|\alpha_i|}{\nu t} 
   \,\exp\Bigl(-\frac{|x-x_i|^2}{5\nu t}\Bigr)~,   
\end{equation}
for all $i \in \{1,\dots,N\}$, all $x \in \real^2$, and all $t > 0$. 
\end{lemma}

The proof of Lemma~\ref{lemdec}, which is borrowed from
\cite{GW05,GG05}, will be reproduced in Section~\ref{app} for the
reader's convenience. For the time being, we just point out that
estimate \eqref{omegaibdd} gives a precise information on the 
$N$-vortex solution in the limit $t \to 0$ for any fixed $\nu > 0$,
but cannot be used to control $\omega_i^\nu(x,t)$ in the limit 
$\nu \to 0$ for fixed $t > 0$, because the constant $K_0(\nu,|\alpha|)$ 
blows up rapidly as $\nu \to 0$.

\medskip
In the very particular case where $N = 1$, we know from \cite{GW05,GGL05}
that $\omega^\nu(x,t)$ is just a suitable translate of Oseen's vortex 
\eqref{Oseen}. From now on, we assume that $N \ge 2$, and we suppose
that the point vortex system \eqref{PW} is well-posed on the time 
interval $[0,T]$. We denote
\begin{equation}\label{dmin}
  d \,=\, \min_{t \in [0,T]}\,\min_{i\neq j}\,|z_i(t) - z_j(t)|
  \,>\, 0~, 
\end{equation}
and we also introduce the turnover time 
\begin{equation}\label{T0def}
  T_0 \,=\, \frac{d^2}{|\alpha|}~, \quad \hbox{where} \quad
  |\alpha| \,=\, |\alpha_1| + \dots + |\alpha_N|~.
\end{equation}
Our goal is to show that, if $\nu > 0$ is sufficiently small, the
$N$-vortex solution $\omega^\nu(x,t)$ looks like a superposition of
$N$ Oseen vortices located at some points $z_1^\nu(t),\dots,
z_N^\nu(t)$ which satisfy the following {\em viscous regularization} 
of the Helmholtz-Kirchhoff system:
\begin{equation}\label{PW2}
  \frac{\D}{\D t} z_i^\nu(t) \,=\, \sum_{j=1}^N \frac{\alpha_j}
  {\sqrt{\nu t}}\, v^G\Bigl(\frac{z_i^\nu(t) - z_j^\nu(t)}
  {\sqrt{\nu t}}\Bigr)~, \qquad z_i^\nu(0) \,=\, x_i~,
\end{equation}
where $v^G$ is given by \eqref{Gdef}. The reason for using that system
instead of \eqref{PW} will be explained at the beginning of
Section~\ref{sec3}. For the moment, we just observe that system
\eqref{PW2} is globally well-posed for positive times, and that the
solutions $z_i^\nu(t)$ are exponentially close to the solutions
$z_i(t)$ of \eqref{PW} if the viscosity $\nu$ is sufficiently small.

\begin{lemma}\label{PWapprox}
Assuming pairwise distinct initial positions $x_1,\dots,x_N$, system 
\eqref{PW2} is globally well-posed for positive times, for any 
value of $\nu > 0$. Moreover, if the solution of \eqref{PW} 
satisfies \eqref{dmin}, there exists $K_1 > 0$ (depending 
on the ratio $T/T_0$) such that
\begin{equation}\label{PWcomp}
  \frac{1}{d}\,\max_{i=1,\dots,N} |z_i^\nu(t) - z_i(t)| \,\le\, 
  K_1 \,\exp\Bigl(-\frac{d^2}{5\nu t}\Bigr)~, 
  \quad \hbox{for all } t \in (0,T]~.
\end{equation}
\end{lemma}

\noindent
The proof of Lemma~\ref{PWapprox} is also postponed to Section~\ref{app}. 
We refer to \cite{NSUW} for a discussion of the validity of system
\eqref{PW2} as a model for the dynamics of interacting vortices.

\medskip To obtain a precise description of each vortex patch
$\omega_i^\nu(x,t)$ in a neighborhood of $z_i^\nu(t)$, we introduce, 
for each $i \in \{1,\dots,N\}$, the self-similar variable
\[
  \xi \,=\, \frac{x - z_i^\nu(t)}{\sqrt{\nu t}}~, 
\]
and we define rescaled functions $w_i^\nu(\xi,t) \in \real$ and 
$v_i^\nu(\xi,t) \in \real^2$ by setting
\begin{equation}\label{wvdef}
  \omega_i^\nu(x,t) \,=\, \frac{\alpha_i}{\nu t}\,w_i^\nu
  \Bigl(\frac{x - z_i^\nu(t)}{\sqrt{\nu t}}\,,\,t\Bigr)~, \qquad
  u_i^\nu(x,t) \,=\, \frac{\alpha_i}{\sqrt{\nu t}}\,v_i^\nu
  \Bigl(\frac{x - z_i^\nu(t)}{\sqrt{\nu t}}\,,\,t\Bigr)~.
\end{equation}
Given a small $\beta \in (0,1)$, which will be specified later, 
we also introduce a weighted $L^2$ space $X$ defined by the 
following norm:
\begin{equation}\label{Xdef}
  \|w\|_X \,=\, \Bigl(\inttwo |w(\xi)|^2 \,e^{\beta |\xi|/4}
  \d\xi\Bigr)^{1/2}~.
\end{equation}
We already know from \cite{GG05} that $w_i^\nu(\xi,t)$ converges to 
the Gaussian profile $G(\xi)$ as $t \to 0$ for any fixed $\nu > 0$.
Our first main result shows that a similar result holds in the 
vanishing viscosity limit, uniformly in time on the interval 
$(0,T]$. 

\begin{theorem}\label{thm3}
Assume that the point vortex system \eqref{PW} is well-posed on 
the time interval $[0,T]$, and let $\omega^\nu(x,t)$ be the 
solution of \eqref{V} with initial data \eqref{mudef}. If 
$\omega^\nu(x,t)$ is decomposed as in \eqref{omdecomp}, then  
the rescaled profiles $w_i^\nu(\xi,t)$ defined by \eqref{wvdef}
satisfy
\begin{equation}\label{strongconv}
  \max_{i=1,\dots,N} \|w_i^\nu(\cdot,t) - G\|_X \,=\ 
  \cO\Bigl(\frac{\nu t}{d^2}\Bigr)~, \quad \hbox{as } \nu \to 0~,
\end{equation}
uniformly for $t \in (0,T]$, where $d$ is given by \eqref{dmin}. 
\end{theorem}

More precisely, the proof shows that there exist positive constants 
$\beta$ and $K_2$, depending on the ratio $T/T_0$, such that
\begin{equation}\label{sconv}
  \max_{i=1,\dots,N} \|w_i^\nu(\cdot,t) - G\|_X \,\le\ 
  K_2\,\frac{\nu t}{d^2}~, 
\end{equation}
for all $t \in (0,T]$, provided $\nu$ is small enough so that $K_2
(\nu T/d^2) \le 1$. This result means that, when the viscosity $\nu$
is small, the $N$-vortex solution $\omega^\nu(x,t)$ looks like a
superposition of $N$ Oseen vortices located at the points
$z_1^\nu(t),\dots, z_N^\nu(t)$, which evolve in time according to
\eqref{PW2}. Since $X \hookrightarrow L^1(\real^2)$ and since
$z_i^\nu(t) \to z_i(t)$ as $\nu \to 0$ by Lemma~\ref{PWapprox},
estimate \eqref{strongconv} implies in particular \eqref{weakconv},
and Theorem~\ref{thm2} is thus a direct consequence of
Theorem~\ref{thm3}.

\medskip Theorem~\ref{thm3} is already very satisfactory, but it does
not seem possible to prove it directly without computing a higher
order approximation of the $N$-vortex solution. This rather surprising
claim will be justified in Section~\ref{sec3}, but for the moment it
can be roughly explained as follows. As is clear from \eqref{Oseen},
the velocity field of Oseen's vortex is very large near the center if
the viscosity $\nu$ is small, with a maximal angular speed of the
order of $|\alpha|/(\nu t)$. As long as the vortex stays isolated, it
does not feel at all the effect of its own velocity field, because of
radial symmetry. However, if a vortex is advected by a non-homogenous
external field, which in our case is produced by the other $N{-}1$
vortices, it will get deformed and, consequently, will start feeling
the influence of its own velocity field. If the ratio $|\alpha|/\nu$
is large, this self-interaction will have a very strong effect, even
if the deformation is quite small. In particular, one may fear that
the vortex gets further deformed, increasing in turn the
self-interaction itself, and that the whole process results in a
violent instability. In fact, this catastrophic scenario does not
happen. Remarkably enough, a rapidly rotating Oseen vortex in an
external field adapts its shape in such a way that the
self-interaction {\em counterbalances} the strain of the external
field \cite{TT65,TK91}. This fundamental observation will be the basis
for our analysis in Section~\ref{sec3}. It explains why one can
observe, in turbulent two-dimensional flows, stable asymmetric
vortices which in a first approximation are simply advected by the
main stream. The same mechanism accounts for the existence of stable
asymmetric Burgers vortices, which are stationary solutions of the
three-dimensional Navier-Stokes equations in a linear strain field
\cite{MKO94,GW07,Ma2,Ma3}.

To compute the self-interactions of the vortices, which play a 
crucial role in the convergence proof, we construct a higher 
order approximation of the $N$-vortex solution in the following
way. For each $i \in \{1,\dots,N\}$ and all $t \in [0,T]$, 
we denote
\begin{equation}\label{wapp}
  w_i^\app(\xi,t) \,=\, G(\xi) + \Bigl(\frac{\nu t}{d^2}\Bigr)
  \Bigl\{\bar F_i(\xi,t) + F_i^\nu(\xi,t)\Bigr\}~, \qquad 
  \xi \in \real^2~,
\end{equation}
where $\bar F_i(\xi,t)$ is a radially symmetric function of 
$\xi$, whose precise expression is given in \eqref{Fibarexp} 
below, and $F_i^\nu(\xi,t)$ is a nonsymmetric correction which 
satisfies
\begin{equation}\label{Finuexp}
  F_i^\nu(\xi,t) \,=\, \frac{d^2}{4\pi}\,\omega(|\xi|) 
  \sum_{j\neq i} \frac{\alpha_j}{\alpha_i}\,\frac{1}{|z_{ij}(t)|^2}
  \Bigl(2\frac{|\xi \cdot z_{ij}(t)|^2}{|\xi|^2 |z_{ij}(t)|^2}
  - 1\Bigr) + \cO\Bigl(\frac{\nu}{|\alpha|}\Bigr)~, 
\end{equation}
where $z_{ij}(t) = z_i^\nu(t) - z_j^\nu(t)$. Here $\omega : (0,\infty)
\to \real$ is a smooth, positive function satisfying $\omega(r)
\approx C_1 r^2$ as $r \to 0$ and $\omega(r) \approx C_2 r^4
e^{-r^2/4}$ as $r \to \infty$ for some $C_1, C_2 > 0$, see 
Eq.~\eqref{Fi0def} below. The right-hand side of \eqref{wapp} is the
beginning of an asymptotic expansion of the rescaled vortex patch
$w_i^\nu(\xi,t)$ in powers of the non-dimensional parameter $(\nu
t)/d^2$. Each term in this expansion can in turn be developped in
powers of $\nu/|\alpha|$. The most important physical effect is due to
the nonsymmetric term $F_i^\nu(\xi,t)$, which describes to leading
order the deformation of the $i^{\rm th}$ vortex due to the influence
of the other vortices. Keeping only that term and using polar
coordinates $\xi = (r\cos\theta,r\sin\theta)$, we can rewrite
\eqref{wapp} in the following simplified form
\begin{equation}\label{wapp2}
  w_i^\app(\xi,t) \,=\, g(r) + \frac{\omega(r)}{4\pi} \sum_{j\neq i} 
  \frac{\alpha_j}{\alpha_i}\,\frac{\nu t}{|z_{ij}(t)|^2}
  \cos\Bigl(2(\theta - \theta_{ij}(t))\Bigr) \,+\, \dots~,
\end{equation}
where $g(|\xi|) = G(\xi)$ and $\theta_{ij}(t)$ is the argument of the 
planar vector $z_{ij}(t) = z_i^\nu(t) - z_j^\nu(t)$. This formula 
allows to compute the principal axes and the eccentricities of 
the vorticity contours, which are elliptical at this level of 
approximation. 

\medskip
Using these notations, our final result can now be stated as follows:

\begin{theorem}\label{thm4}
Assume that the point vortex system \eqref{PW} is well-posed on 
the time interval $[0,T]$, and let $\omega^\nu(x,t)$ be the 
solution of \eqref{V} with initial data \eqref{mudef}. If 
$\omega^\nu(x,t)$ is decomposed as in \eqref{omdecomp}, then  
the rescaled profiles $w_i^\nu(\xi,t)$ defined by \eqref{wvdef}
satisfy
\begin{equation}\label{strongconv2}
  \max_{i=1,\dots,N} \|w_i^\nu(\cdot,t) - w_i^\app(\cdot,t)\|_X 
  \,=\, \cO\Bigl(\Bigl(\frac{\nu t}{d^2}\Bigr)^{3/2}\Bigr)~, 
  \quad \hbox{as } \nu \to 0~,
\end{equation}
uniformly for $t \in (0,T]$, where $w_i^\app$ is given by 
\eqref{wapp}. 
\end{theorem}

As is clear from \eqref{strongconv}, \eqref{wapp},
\eqref{strongconv2}, Theorem~\ref{thm4} implies immediately
Theorem~\ref{thm3}, hence also Theorem~\ref{thm2}. Note that the
convergence result is now accurate enough so that the first order
corrections to the Gaussian profile in \eqref{wapp} are much larger,
for small viscosities, than the remainder terms which are summarized
in the right-hand side of \eqref{strongconv2}. This means in
particular that the deformations of the interacting vortices are
really given, to leading order, by \eqref{wapp2}. According to that
formula, each vortex adapts its shape {\em instantaneously} to the
relative positions of the other vortices, without oscillations or
inertia. In this sense, point vortices can be considered as
well-prepared initial data, and the $N$-vortex solution is an example
of the ``metastable regime'' described in the introduction. In
contrast, one should mention that the first order radially symmetric
corrections $\bar F_i(\xi,t)$ do not only depend on the instantaneous
vortex positions, but on the whole history of the system, see
\eqref{Fibarexp}.

The rest of this paper is devoted to the proof of Theorem~\ref{thm4},
which is divided in two main steps. In Section~\ref{sec3}, we
construct an approximate solution of our system, with the property
that the associated residuum is extremely small in the vanishing
viscosity limit. This approximation differs from \eqref{wapp} by
higher order corrections which are necessary to reach the desired
accuracy, but will eventually be absorbed in the right-hand side of
\eqref{strongconv2}. As is explained in Section~\ref{3.3}, a
difficulty in this construction comes from the fact that the radially
symmetric and the nonsymmetric terms in the approximate solution
$w_i^\app(\xi,t)$ have a different origin, and play a different role.
Once a suitable approximation has been constructed, our final task is
to control the remainder $w_i(\xi,t) - w_i^\app(\xi,t)$ uniformly in
time and in the viscosity parameter $\nu$. This will be done in
Section~\ref{sec4}, using appropriate energy estimates.  A major
technical difficulty comes here from the fact that we do not want to
assume that the ratio $T/T_0$ is small. If we did so, the proof would
be considerably simpler, and we could replace the weight
$e^{\beta|\xi|/4}$ by $e^{\beta|\xi|^2/4}$ in the definition
\eqref{Xdef} of our function space $X$, thus improving our convergence
result. In the general case, however, we have to use a rather delicate
energy estimate involving time-dependent weights $p_i(\xi,t)$, which
will be constructed in Section~\ref{4.1}. The price to pay is a
slightly weaker control of the remainder, and the fact that all our
constants, such as $K_2$ and $\beta$ in \eqref{sconv}, have a bad
dependence on $T$ if $T \gg T_0$.

In conclusion, our results show that the vanishing viscosity limit can
be rigorously controlled in the particular case of point vortices, due
to the remarkable dynamic and strucural stability properties of the
Oseen vortices. These properties, which were established in 
\cite{GW05,GW07,Ma1}, play a crucial role both in the construction
of the approximate solution in Section~\ref{sec3}, and in the 
energy estimates of Section~\ref{sec4}.

\section{Construction of an approximate solution}\label{sec3}

In this section, we show how to construct an asymptotic expansion
of the $N$-vortex solution in the vanishing viscosity limit. 
Our starting point is the evolution system satisfied by the 
rescaled profiles $w_i^\nu(\xi,t)$, $v_i^\nu(\xi,t)$, which 
from now on will be denoted by $w_i(\xi,t)$, $v_i(\xi,t)$ 
for simplicity. Replacing \eqref{wvdef} into \eqref{convecdiff}, 
we obtain, for $i = 1,\dots,N$, 
\begin{equation}\label{wieq}
  t\partial_t w_i(\xi,t) + \left\{\sum_{j=1}^N \frac{\alpha_j}{\nu}
  \,v_j\Bigl(\xi + \frac{z_{ij}(t)}{\sqrt{\nu t}}\,,\,t\Bigr)
  - \sqrt{\frac{t}{\nu}}\,z_i'(t)\right\}\cdot \nabla w_i(\xi,t) 
  \,=\, (\cL w_i)(\xi,t)~,
\end{equation}
where 
\begin{equation}\label{cLdef}
  \cL w \,=\, \Delta w + \frac12 \xi\cdot \nabla w + w~.
\end{equation}
Here and in what follows, we denote the vortex positions by 
$z_i(t)$ instead of $z_i^\nu(t)$, to keep the formulas as simple 
as possible. We also recall that $z_{ij}(t) = z_i(t)-z_j(t)$. 

The initial value problem for system~\eqref{wieq} at time $t = 0$
is not well-posed, because the time derivative appears in 
the singular form $t\partial_t$. A convenient way to avoid this
difficulty is to introduce a new variable $\tau = \log(t/T)$, 
so that $\partial_\tau = t\partial_t$. With this parametrization, 
the solution of \eqref{wieq} given by Lemma~\ref{lemdec} is defined
for all $\tau \in (-\infty,0]$, and converges to the profile $G$ of 
Oseen's vortex as $\tau \to -\infty$, see \cite[Proposition~4.5]{GG05}. 
For simplicity, we keep here the original time $t$, because 
this is the natural variable for the ODE system~\eqref{PW}
or \eqref{PW2}. 

\subsection{Residuum of the naive approximation}\label{3.1}

If $w(\xi,t) = (w_1(\xi,t),\dots,w_N(\xi,t))$ is an approximate 
solution of system \eqref{wieq}, we define the {\em residuum} 
$R(\xi,t) = (R_1(\xi,t),\dots,R_N(\xi,t))$ of this approximation 
by
\[
  R_i(\xi,t) \,=\, 
  t\partial_t w_i(\xi,t) + \biggl\{\sum_{j=1}^N \frac{\alpha_j}{\nu}
  \,v_j\Bigl(\xi + \frac{z_{ij}(t)}{\sqrt{\nu t}}\,,\,t\Bigr)
  - \sqrt{\frac{t}{\nu}}\,z_i'(t)\biggr\}\cdot \nabla w_i(\xi,t) 
  - (\cL w_i)(\xi,t)~,
\]
for all $i \in \{1,\dots,N\}$. Here and in the sequel, it is always 
understood that $v_i(\xi,t)$ is the velocity field corresponding to 
$w_i(\xi,t)$ via the Biot-Savart law \eqref{BS}. 

In view of Theorem~\ref{thm3}, the solution of \eqref{wieq} we are 
interested in satisfies $w_i(\xi,t) \approx G(\xi)$ and $v_i(\xi,t) 
\approx v^G(\xi)$ for all $t \in (0,T]$, if $\nu$ is sufficiently
small. Since $\partial_t G = \cL G = 0$, the residuum of this naive 
approximation (where $w_i(\xi,t) = G(\xi)$ for all $i \in \{1,\dots,N\}$)
is
\begin{equation}\label{R0def}
  R_i^{(0)}(\xi,t) \,=\, \biggl\{\sum_{j=1}^N \frac{\alpha_j}{\nu}
  \,v^G\Bigl(\xi + \frac{z_{ij}(t)}{\sqrt{\nu t}}\Bigr)
  - \sqrt{\frac{t}{\nu}}\,z_i'(t)\biggr\}\cdot \nabla G(\xi)~.
\end{equation}
This expression looks singular in the limit $\nu \to 0$, but
the problem can be eliminated by an appropriate choice of the 
vortex positions $z_1(t),\dots,z_N(t)$. Indeed, in \eqref{R0def}, 
the quantity inside the curly brackets $\{\cdot\}$ vanishes for 
$\xi = 0$ if we set
\begin{equation}\label{PW3}
  z_i'(t) \,=\, \sum_{j=1}^N \frac{\alpha_j}{\sqrt{\nu t}}
  \, v^G\Bigl(\frac{z_{ij}(t)}{\sqrt{\nu t}}\Bigr)~, \qquad
  i = 1,\dots,N~.
\end{equation}
This is exactly the regularized point vortex system, which was already
introduced in \eqref{PW2}. 

From now on we always assume, as in Theorems~\ref{thm3} and
\ref{thm4}, that the original point vortex system \eqref{PW} is
well-posed on the time interval $[0,T]$. For any $\nu > 0$, we denote
by $z_1(t),\dots,z_N(t)$ the solution of \eqref{PW3} with initial data
$x_1,\dots,x_N$. In view of Lemma~\ref{PWapprox}, we can assume that
this solution satisfies \eqref{dmin} for some $d > 0$ (independent of
$\nu$), provided $\nu$ is sufficiently small. Replacing \eqref{PW3}
into \eqref{R0def}, we obtain the following expression of the residuum
\begin{equation}\label{R0def2}
  R_i^{(0)}(\xi,t) \,=\, \sum_{j\neq i} \frac{\alpha_j}{\nu}
  \biggl\{v^G\Bigl(\xi + \frac{z_{ij}(t)}{\sqrt{\nu t}}\Bigr)
    - v^G\Bigl(\frac{z_{ij}(t)}{\sqrt{\nu t}}\Bigr)\biggr\}\cdot 
  \nabla G(\xi)~. 
\end{equation}
Remark that the sum now runs on the indices $j \neq i$, because 
the term corresponding to $j = i$ vanishes. Our first task is 
to compute an asymptotic expansion of the right-hand side of 
\eqref{R0def2} in the vanishing viscosity limit. 

\begin{proposition}\label{R0expand}
For $i = 1,\dots,N$, we have
\begin{equation}\label{R0exp}
  R_i^{(0)}(\xi,t) \,=\, \frac{\alpha_i t}{d^2}\biggl\{
  A_i(\xi,t) + \Bigl(\frac{\nu t}{d^2}\Bigr)^{1/2}
  B_i(\xi,t) +  \Bigl(\frac{\nu t}{d^2}\Bigr) C_i(\xi,t) 
  + \tilde R_i^{(0)}(\xi,t)\biggr\}~,
\end{equation}
for all $\xi \in \real^2$ and all $t \in (0,T]$, where
\begin{eqnarray}\nonumber
  A_i(\xi,t) &=& \frac{d^2}{2\pi} \sum_{j \neq i}
  \frac{\alpha_j}{\alpha_i}\,\frac{(\xi\cdot z_{ij}(t))
  (\xi\cdot z_{ij}(t)^\perp)}{|z_{ij}(t)|^4}\,G(\xi)~,\\ \label{ABC}
  B_i(\xi,t) &=& \frac{d^3}{4\pi} \sum_{j \neq i}
  \frac{\alpha_j}{\alpha_i}\,\frac{(\xi\cdot z_{ij}(t)^\perp)}
  {|z_{ij}(t)|^6}\,\Bigl(|\xi|^2 |z_{ij}(t)|^2 - 4 
  (\xi\cdot z_{ij}(t))^2\Bigr)\,G(\xi)~,\\ \nonumber
  C_i(\xi,t) &=& \frac{d^4}{\pi} \sum_{j \neq i}
  \frac{\alpha_j}{\alpha_i}\,\frac{(\xi\cdot z_{ij}(t))
  (\xi\cdot z_{ij}(t)^\perp)}{|z_{ij}(t)|^8}\Bigl(
  2(\xi\cdot z_{ij}(t))^2 - |\xi|^2 |z_{ij}(t)|^2\Bigr)\,G(\xi)~.
\end{eqnarray}
Moreover, for any $\gamma < 1$, there exists $C > 0$ such that
\begin{equation}\label{R0rem}
   |\tilde R_i^{(0)}(\xi,t)| \,\le\, C\Bigl(\frac{\nu t}{d^2}
   \Bigr)^{3/2}\,e^{-\gamma|\xi|^2/4}~, \qquad \xi \in \real^2~, 
   \quad 0 < t \le T~. 
\end{equation}
\end{proposition}

\noindent{\bf Remarks.} Proposition~\ref{R0expand} provides an
expansion of the residuum $R_i^{(0)}(\xi,t)$ in powers of the
dimensionless parameter $(\nu t)^{1/2}/d$, where $(\nu t)^{1/2}$ is a
diffusion length which gives the typical size of each vortex at time
$t$, and $d$ is the minimal distance between the vortex centers, see
\eqref{dmin}. As is clear from \eqref{ABC}, the parameter $d$ in
\eqref{R0exp} has been introduced rather artificially, to ensure that
the quantities $A_i$, $B_i$, $C_i$ are dimensionless, as is the
residuum itself. The proof will show that an expansion of the form
\eqref{R0exp} can be performed to arbitrarily high orders, but for
simplicity we keep only the terms which will be necessary to prove
Theorems~\ref{thm3} and \ref{thm4}.  Finally, we remark that the
prefactor $\alpha_i t/d^2$ in \eqref{R0exp} is bounded by $t/T_0$,
where $T_0$ is the turnover time introduced in \eqref{T0def}.

\medskip\proof Fix $\gamma \in (0,1)$, and let $\gamma_1 = 
1 - \gamma$. Since $\nabla G(\xi) = -\frac12 \xi G(\xi)$ 
decreases rapidly as $|\xi| \to \infty$, it is clear that 
the residuum \eqref{R0def2} is extremely small if $|\xi|$ 
is large. For instance, if $|\xi| \ge d/(2\sqrt{\nu t})$, 
we can bound
\begin{align*}
  |R_i^{(0)}(\xi,t)| \,&\le\, \frac{|\alpha|}{\nu}\,\|v^G\|_{L^\infty}
  |\xi|G(\xi) \,=\, \frac{1}{4\pi}\,\frac{|\alpha|}{\nu} 
  \,\|v^G\|_{L^\infty}|\xi|\,e^{-\gamma_1|\xi|^2/4}\,e^{-\gamma|\xi|^2/4}\\
  \,&\le\, C\,\frac{|\alpha|t}{d^2}\,\Bigl(\frac{d^2}{\nu t}
  \Bigr)^{3/2}\exp\Bigl(-\frac{\gamma_1 d^2}{16\nu t}\Bigr)
  \,e^{-\gamma|\xi|^2/4}~,
\end{align*}
where in the last inequality we have used the fact that 
$|\xi|\,e^{-\gamma_1|\xi|^2/4}$ is a decreasing function of
$|\xi|$ when $|\xi| \gg 1$. A similar argument shows that 
\[
  |A_i(\xi,t)| + |B_i(\xi,t)| + |C_i(\xi,t)| \,\le\, 
  C\,\Bigl(\frac{d^2}{\nu t}\Bigr)^2\exp\Bigl(-
  \frac{\gamma_1 d^2}{16\nu t}\Bigr)\,e^{-\gamma|\xi|^2/4}~,
\]
if $|\xi| \ge d/(2\sqrt{\nu t})$. We conclude that expansion 
\eqref{R0exp} holds in that region, with a remainder term 
satisfying a much better estimate than \eqref{R0rem}. 

We now consider the case where $|\xi| \le d/(2\sqrt{\nu t})$. Since 
$|z_{ij}(t)| = |z_i(t)-z_j(t)| \ge d$ by \eqref{dmin} if $i \neq j$, 
we have 
\begin{equation}\label{auxbounds}
  \Bigl|\frac{z_{ij}(t)}{\sqrt{\nu t}}\Bigr| \,\ge\, 
  \frac{d}{\sqrt{\nu t}} \,\ge\, 2|\xi|~, 
  \quad \hbox{and}\quad 
  \Bigl|\xi + \frac{z_{ij}(t)}{\sqrt{\nu t}}\Bigr| \,\ge\, 
  \frac{d}{2\sqrt{\nu t}} \,\ge\, |\xi|~.
\end{equation}
To estimate the right-hand side of \eqref{R0def2}, we have to compute
the difference $v^G(\xi+\eta) - v^G(\eta)$, with $\eta = 
z_{ij}(t)/\sqrt{\nu t}$. Using definition \eqref{Gdef}, we obtain 
the identity 
\begin{equation}\label{V1V2}
  v^G(\xi+\eta) - v^G(\eta) \,=\, \frac{1}{2\pi}\Bigl(V_1(\xi,\eta) + 
  V_2(\xi,\eta)\Bigr)~, \qquad \xi,\eta \in \real^2~,
\end{equation}
where
\[
  V_1(\xi,\eta) \,=\, \frac{(\xi+\eta)^\perp}{|\xi+\eta|^2} 
  - \frac{\eta^\perp}{|\eta|^2}~, \qquad
  V_2(\xi,\eta) \,=\, \frac{\eta^\perp}{|\eta|^2}
  \,e^{-|\eta|^2/4} - \frac{(\xi+\eta)^\perp}{|\xi+\eta|^2}
  \,e^{-|\xi+\eta|^2/4}~.
\]
In particular, it follows from \eqref{auxbounds} that
\begin{equation}\label{V2exp}
  \Bigl|V_2\Bigl(\xi,\frac{z_{ij}(t)}{\sqrt{\nu t}}\Bigr)\Bigr|
  \,\le\, C\,\Bigl(\frac{\nu t}{d^2}\Bigr)^{1/2}\exp\Bigl(-
  \frac{d^2}{16\nu t}\Bigr)~,  
\end{equation}
hence the contributions of the term $V_2$ to the residuum 
\eqref{R0def2} are exponentially small and can be incorporated
in the remainder term. To compute $V_1$, we use the following 
elementary lemma: 

\begin{lemma}\label{Rexpand}
For all $\xi,\eta \in \real^2$ with $|\xi| < |\eta|$, we have
\begin{equation}\label{V11}
  \xi\cdot V_1(\xi,\eta) \,=\, \sum_{n=2}^\infty (-1)^{n-1}
  \,\frac{|\xi|^n}{|\eta|^n}\,\sin(n(\theta-\phi))~,
\end{equation} 
where $\theta$ denotes the polar argument of $\xi$ and $\phi$ 
the argument of $\eta$. 
\end{lemma}

For completeness, the proof of Lemma~\ref{Rexpand} will be given in
Section~\ref{app}. Applying \eqref{V11} with $\eta =
z_{ij}(t)/\sqrt{\nu t}$ and keeping only the first three terms in the
expansion, we obtain
\begin{align}\nonumber
  \xi\cdot V_1\Bigl(\xi,\frac{z_{ij}}{\sqrt{\nu t}}\Bigr)
  \,=\, &-\frac{|\xi|^2 \nu t}{|z_{ij}|^2}\,\sin(2(\theta-\phi))
  + \frac{|\xi|^3 (\nu t)^{3/2}}{|z_{ij}|^3}\,\sin(3(\theta-\phi)) \\
  \label{xiV}
  & -\frac{|\xi|^4 (\nu t)^2}{|z_{ij}|^4}\,\sin(4(\theta-\phi))
  +  \cO\Bigl(\frac{|\xi|^5 (\nu t)^{5/2}}{|z_{ij}|^5}\Bigr)~,
\end{align}
where $\theta-\phi$ is the signed angle between $\xi$ and $z_{ij}$.
In particular, we have the relations
\[
  \sin(\theta-\phi) \,=\, \frac{\xi\cdot z_{ij}^\perp}{|\xi||z_{ij}|}~, 
  \quad \cos(\theta-\phi) \,=\, \frac{\xi\cdot z_{ij}}{|\xi||z_{ij}|}~,  
\]
in terms of which the higher order trigonometric expressions 
appearing in \eqref{xiV} can be computed using the well-known formulas
$\sin(2\alpha) = 2\sin(\alpha) \cos(\alpha)$, $\sin(3\alpha) =
\sin(\alpha)(4\cos^2(\alpha)-1)$, and $\sin(4\alpha) =
4\sin(\alpha)\cos(\alpha)(2\cos^2(\alpha)-1)$. Summarizing, we
have shown that
\begin{align*}
  \sum_{j\neq i}\frac{\alpha_j}{\nu}\,\frac{1}{2\pi} &V_1\Bigl(\xi,
  \frac{z_{ij}(t)}{\sqrt{\nu t}}\Bigr)\cdot \nabla G(\xi) 
  \,=\, -\frac{1}{4\pi}\sum_{j\neq i}\frac{\alpha_j}{\nu}
  \,\xi\cdot V_1\Bigl(\xi,\frac{z_{ij}(t)}{\sqrt{\nu t}}\Bigr) G(\xi) \\
  \,&=\, \frac{\alpha_i t}{d^2}\biggl\{
  A_i(\xi,t) + \Bigl(\frac{\nu t}{d^2}\Bigr)^{1/2}
  B_i(\xi,t) +  \Bigl(\frac{\nu t}{d^2}\Bigr) C_i(\xi,t) 
  + \hat R_i(\xi,t)\biggr\}~,
\end{align*}  
for $|\xi| \le d/(2\sqrt{\nu t})$, where $A_i$, $B_i$, $C_i$ are given
by \eqref{ABC} and $\hat R_i(\xi,t)$ satisfies the bound
\eqref{R0rem}. As was already observed, the same result holds if we
replace $V_1$ with $V_1 + V_2$, so we conclude that expansion
\eqref{R0exp} is valid in the region $|\xi| \le d/(2\sqrt{\nu t})$
too. The proof of Proposition~\ref{R0expand} is thus complete. \QED

It is important to remark that the residuum $R_i^{(0)}(\xi,t)$ does
not converge to zero as $\nu \to 0$, because of the leading order term
$A_i(\xi,t)$. If we decompose the solution of \eqref{wieq} as
$w_i(\xi,t) = G(\xi) + \tilde w_i(\xi,t)$, the equation for $\tilde
w_i(\xi,t)$ will contain a source term of size $\cO(1)$ as $\nu \to
0$, and we therefore expect that the remainder $\tilde w_i(\xi,t)$
itself will be of size $\cO(1)$ after a short time.  But, as is easily
verified, the equation for $\tilde w_i(\xi,t)$ contains nonlinear
terms with a prefactor of size $\cO(\nu^{-1})$, and such terms cannot
be controlled in the vanishing viscosity limit if $\tilde w_i(\xi,t)$
is $\cO(1)$. This is the reason why it is necessary to construct a more
precise approximate solution of \eqref{wieq}, with a sufficiently
small residuum, in order to desingularize the equation for the
remainder.

\subsection{First order approximation}\label{3.2}

We look for an approximate solution of \eqref{wieq} of the form
\begin{equation}\label{wvfirst}
  w_i^\app(\xi,t) \,=\, G(\xi) + \Bigl(\frac{\nu t}{d^2}\Bigr)
  F_i(\xi,t)~, \qquad 
  v_i^\app(\xi,t) \,=\,  v^G(\xi) + \Bigl(\frac{\nu t}{d^2}\Bigr)
  v^{F_i}(\xi,t)~,
\end{equation} 
where, for each $i \in \{1,\dots,N\}$, $F_i(\xi,t)$ is a smooth 
vorticity profile to be determined, and $v^{F_i}(\xi,t)$ is the 
velocity field obtained from $F_i(\xi,t)$ via the Biot-Savart 
law \eqref{BS}. In fact, we shall need later a more precise 
approximation of the $N$-vortex solution, but we prefer starting 
with \eqref{wvfirst} to describe the procedure in a relatively 
simple setting. Our goal is to chose the profile $F_i(\xi,t)$ 
so as to minimize the residuum of our approximation, which is
\[
  R_i^{(1)}(\xi,t) \,=\, (t\partial_t -\cL) w_i^\app(\xi,t)  
  + \sum_{j=1}^N \frac{\alpha_j}{\nu}\biggl\{v_j^\app\Bigl(\xi + 
  \frac{z_{ij}(t)}{\sqrt{\nu t}}\,,t\Bigr) - v^G
  \Bigl(\frac{z_{ij}(t)}{\sqrt{\nu t}}\Bigr)\biggr\}\cdot 
  \nabla w_i^\app(\xi,t)~.
\]
Using \eqref{wvfirst} and the definition \eqref{R0def2} of 
$R_i^{(0)}(\xi,t)$, we find after some calculations
\begin{equation}\label{R1def}
  R_i^{(1)}(\xi,t) \,=\, R_i^{(0)}(\xi,t) + \frac{\alpha_i t}{d^2}
  \Bigl(v^G \cdot \nabla F_i + v^{F_i}\cdot\nabla G\Bigr)(\xi,t)
  + \tilde R_i^{(1)}(\xi,t)~, 
\end{equation}
where
\begin{align*}
  \tilde R_i^{(1)}(\xi,t) \,&=\, \Bigl(\frac{\nu t}{d^2}\Bigr)
  \Bigl(t\partial_t F_i + F_i -\cL F_i\Bigr)(\xi,t) + 
  \Bigl(\frac{\nu t}{d^2}\Bigr)\sum_{j=1}^N \frac{\alpha_j t}{d^2}
  \,v^{F_j}\Bigl(\xi + \frac{z_{ij}}{\sqrt{\nu t}}\,,t\Bigr)
  \cdot\nabla F_i(\xi,t)\\   
  \,&+\, \sum_{j\neq i} \frac{\alpha_j t}{d^2}\biggl\{
  \Bigl(v^G\Bigl(\xi + \frac{z_{ij}}{\sqrt{\nu t}}\Bigr) - 
  v^G\Bigl(\frac{z_{ij}}{\sqrt{\nu t}}\Bigr)\Bigr)\cdot\nabla 
  F_i(\xi,t) + v^{F_j}\Bigl(\xi + \frac{z_{ij}}{\sqrt{\nu t}}\,,t\Bigr)   
  \cdot \nabla G(\xi)\biggr\}~.
\end{align*}
It is easy to check, at least formally, that $\tilde R_i^{(1)}(\xi,t)$
is $\cO(\nu t/d^2)$ as $\nu \to 0$ (this calculation will be 
done rigorously later, when the profile $F_i$ will be determined).
So it follows from \eqref{R0exp} and \eqref{R1def} that
\begin{equation}\label{R1exp}
  R_i^{(1)}(\xi,t) \,=\, \frac{\alpha_i t}{d^2}
  \Bigl(A_i + v^G \cdot \nabla F_i + v^{F_i}\cdot\nabla G\Bigr)
  (\xi,t) + \cO\Bigl(\Bigl(\frac{\nu t}{d}\Bigr)^{\frac12}\Bigr)~.
\end{equation}
To minimize the residuum, it is natural to impose $A_i + v^G \cdot \nabla
F_i + v^{F_i}\cdot\nabla G = 0$. As we shall see, this ``elliptic''
equation has a solution, but does not completely determine the profile
$F_i$.

To prove this claim, we first introduce some notation. Let $Y$ denote 
the Hilbert space
\begin{equation}\label{Ydef}
  Y \,=\, \Bigl\{w \in L^2(\real^2) \,\Big|\, \inttwo
  |w(\xi)|^2 \,e^{|\xi|^2/4}\d\xi < \infty\Bigr\}~,
\end{equation}
equipped with the scalar product
\begin{equation}\label{Yprod}
  (w_1,w_2)_Y \,=\, \inttwo w_1(\xi) w_2(\xi) \,e^{|\xi|^2/4}\d\xi~.
\end{equation}
We consider the linear operator $\Lambda : D(\Lambda) \to Y$ defined by 
$D(\Lambda) = \{w \in Y\,|\, v^G \cdot \nabla w \in Y\}$ and
\begin{equation}\label{Lamdef}
  \Lambda w \,=\, v^G\cdot\nabla w + v\cdot \nabla G~,
  \qquad w \in D(\Lambda)~,
\end{equation}
where (as always) $v$ denotes the velocity field obtained from $w$ via
the Biot-Savart law \eqref{BS}. With these notations, the equation we
have to solve becomes $\Lambda F_i + A_i = 0$, and we are therefore
interested in computing the (partial) inverse of $\Lambda$ on
appropriate subspaces. This operator has been extensively studied, 
because it plays a prominent role in the stability properties of the 
Oseen vortices and the construction of asymmetric Burgers vortices 
\cite{GW05,GW07,Ma1,Ma2,Ma3}. In particular, we have

\begin{proposition}\label{Lamprop} \cite{GW05,Ma1} The operator
$\Lambda$ is skew-adjoint in $Y$, so that $\Lambda^* = -\Lambda$. 
Moreover,
\begin{equation}\label{KerLam}
  \Ker(\Lambda) \,=\, Y_0 \oplus \{\beta_1\partial_1 G + \beta_2
  \partial_2 G \,|\, \beta_1, \beta_2 \in \real\}~,
\end{equation}
where $Y_0 \subset Y$ is the subspace of all radially symmetric 
functions. 
\end{proposition}

If $w \in Y_0$, the corresponding velocity field $v$ satisfies 
$\xi\cdot v(\xi) = 0$, hence \eqref{Lamdef} immediately implies 
that $\Lambda w = 0$. On the other hand, if we differentiate the
identity $v^G\cdot\nabla G = 0$ with respect to $\xi_1$ and $\xi_2$, 
we obtain $\Lambda (\partial_1 G) = \Lambda (\partial_2 G) = 0$.
So the right-hand side of \eqref{KerLam} is certainly contained 
in the kernel of $\Lambda$, and the converse inclusion was proved
in \cite{Ma1}. On the other hand, the fact that $\Lambda$ is 
skew-adjoint in $Y$ implies that $\Ker(\Lambda) = \Ran(\Lambda)^\perp$, 
so we know that the range of $\Lambda$ is a dense subspace of
$\Ker(\Lambda)^\perp$. We shall now prove that $\Ran(\Lambda)$ 
contains $\Ker(\Lambda)^\perp \cap Z$, where 
\[
  Z \,=\, \Bigl\{w : \real^2 \to \real\,\Big|\, e^{|\xi|^2/8} w \in 
  \cS(\real^2)\Bigr\} \,\subset\, Y~.
\]
Here $\cS(\real^2)$ denotes the space of all smooth, rapidly decreasing 
functions on $\real^2$. 

As was observed e.g. in \cite{GW05}, the operator $\Lambda$ 
commutes with the group $SO(2)$ of all rotations about the origin.
It is thus natural to decompose 
\[
  Y \,=\, \mathop{\oplus}\limits_{n=0}^\infty Y_n \,=\, 
  \mathop{\oplus}\limits_{n=0}^\infty P_n Y~,
\]
where $P_n$ is the orthogonal projection defined in polar 
coordinates $(r,\theta)$ by the formula
\[
  (P_n w)(r\cos\theta,r\sin\theta) \,=\, \frac{2-\delta_{n,0}}{2\pi}
  \int_0^{2\pi} w(r\cos\theta',r\sin\theta')\cos(n(\theta-\theta'))
  \d\theta'~, \quad n \in \natural~.
\]
Then $Y_0 = P_0 Y$ is the subspace of all radially symmetric 
functions, and for $n \ge 1$ the subspace $Y_n = P_n Y$ contains 
all functions of the form $w(r\cos\theta,r\sin\theta) = 
a_1(r)\cos(n\theta) + a_2(r)\sin(n\theta)$. It is not difficult to 
verify that the projections $P_n$ commute with $\Lambda$ for all 
$n \in \natural$, and explicit formulas for the restrictions 
$\Lambda_n = P_n \Lambda P_n$ are given in \cite[Section~4.1.1]{GW05}.
To formulate the main technical result of this section, we shall use 
the following notations:
\begin{equation}\label{gfh}
  g(r) \,=\, \frac{1}{4\pi}\,e^{-r^2/4}~, \quad 
  \phi(r) \,=\, \frac{1}{2\pi r^2}(1 - e^{-r^2/4})~, \quad 
  h(r) \,=\, \frac{g(r)}{2\phi(r)} \,=\, \frac{r^2/4}{e^{r^2/4}-1}~, 
  \quad r > 0~.
\end{equation}

\begin{lemma}\label{Laminv}
  If $z \in Y_n \cap Z$ for some $n \ge 2$, there exists a unique $w
  \in Y_n \cap Z$ such that $\Lambda w = z$.  In particular, if $z =
  a(r)\sin(n\theta)$, then $w = -\omega(r) \cos(n\theta)$, where
\begin{equation}\label{omexp}
  \omega(r) \,=\, h(r)\Omega(r) + \frac{a(r)}{n\phi(r)}~, 
  \qquad r > 0~,
\end{equation}
and $\Omega : (0,\infty) \to \real$ is the unique solution 
of the differential equation
\begin{equation}\label{Omexp}
  -\frac1r (r\Omega'(r))' + \Bigl(\frac{n^2}{r^2} - h(r)\Bigr)
  \Omega(r) \,=\, \frac{a(r)}{n\phi(r)}~, \qquad r > 0~,  
\end{equation}
such that $\Omega(r) = \cO(r^n)$ as $r \to 0$ and $\Omega(r) = 
\cO(r^{-n})$ as $r \to \infty$. 
\end{lemma}

\medskip\noindent{\bf Remarks.}\\
{\bf 1.} By rotation invariance, if $z = a(r)\cos(n\theta)$, then $w =
\omega(r)\sin(n\theta)$, and the relation between $\omega$ and $a$ is
unchanged. The general case where $z = a_1(r)\cos(n\theta)
+ a_2(r)\sin(n\theta)$ follows by linearity.\\[1pt]
{\bf 2.} The conclusion of Lemma~\ref{Laminv} is wrong for $n = 1$.
Indeed, since $\partial_k G = -\frac12 \xi_kG$ for $k = 1,2$, it is
clear that $\partial_1 G, \partial_2 G \in Y_1 \cap Z$, but
Proposition~\ref{Lamprop} asserts that these functions belong to
$\Ker(\Lambda) = \Ran(\Lambda)^\perp$. However, if $z \in Y_1 \cap Z$
satisfies $(z,\partial_kG)_Y = 0$ for $k = 1,2$, one can show that
there exists a unique $w \in Y_1 \cap Z \cap \Ker(\Lambda)^\perp$ such
that $\Lambda w = z$. This result will not be needed in what follows,
so we omit the proof.\\[1pt]
{\bf 3.} If $w$ is as in Lemma~\ref{Laminv} and if $v$ is the 
velocity field associated to $w$ via the Biot-Savart law 
\eqref{BS}, the proof will show that $v$ is smooth and satisfies
\begin{equation}\label{vnbdd}
  |v(\xi)| \,=\, \cO(|\xi|^{n-1}) \quad \hbox{as } \xi \to 0~, 
  \quad \hbox{and} \quad
  |v(\xi)| \,=\, \cO(|\xi|^{-n-1}) \quad \hbox{as } |\xi| \to \infty~. 
\end{equation}

\medskip\noindent{\bf Proof of Lemma~\ref{Laminv}.}
A particular case of Lemma~\ref{Laminv} was proved in 
\cite[Proposition~3.1]{GW07}. Since the general case is quite 
similar, we just indicate here the main steps and refer to 
\cite{GW07} for further details. 

Assume that $z \in Y_n \cap Z$ for some $n \ge 2$.  By
Proposition~\ref{Lamprop}, we have $z \in \Ker(\Lambda)^\perp =
\overline{\Ran(\Lambda)}$. Our task is to verify that $z \in
\Ran(\Lambda)$, and that there exists a unique $w \in Y_n \cap Z$
such that $\Lambda w = z$. Without loss of generality, we assume 
that $z = a(r)\sin(n\theta)$. Then $a : \real_+ \to \real$ is
a smooth function with the property that $e^{r^2/8}a(r)$ decays
rapidly as $r \to \infty$. Furthermore, we can write $a(r) = r^n
A(r^2)$, where $A : [0,\infty) \to \real$ is again a smooth
function. In particular, we have $a(r) = \cO(r^n)$ as $r \to 0$. We
look for a solution $w$ of the form $w = -\omega(r)\cos(n\theta)$. The
corresponding stream function, which is defined by the relation
$-\Delta \Psi = w$, satisfies $\Psi = -\Omega(r)\cos(n\theta)$, where
$\Omega$ is the unique regular solution of the differential equation
\begin{equation}\label{auxdiff}
  -\frac1r (r\Omega'(r))' + \frac{n^2}{r^2}\,\Omega(r) \,=\, 
  \omega(r)~, \qquad r > 0~.  
\end{equation}
Moreover, the velocity field $v = -\nabla^\perp \Psi$ has the 
following expression
\[
  v \,=\, \frac{n}{r}\,\Omega(r)\sin(n\theta)\,\mathbf{e}_r + 
  \Omega'(r)\cos(n\theta)\,\mathbf{e}_\theta~, 
\]
where $\mathbf{e}_r$, $\mathbf{e}_\theta$ are the unit vectors
in the radial and azimuthal directions, respectively. Thus, 
using definitions \eqref{Gdef} and \eqref{gfh}, we obtain
\[
  v\cdot\nabla G \,=\, -\frac{n}{2}\,\Omega(r)g(r)\sin(n\theta)~, 
  \quad \hbox{and} \quad
  v^G\cdot \nabla w \,=\, n \omega(r)\phi(r) \sin(n\theta)~.
\]
In particular, we see that $\Lambda w = z$ if and only if 
$-\frac{n}{2}\Omega g + n\phi\omega = a$, which is \eqref{omexp}.
Furthermore, combining \eqref{omexp} and \eqref{auxdiff},
we obtain the differential equation \eqref{Omexp} which 
determines $\Omega$. 

It remains to verify that \eqref{Omexp} has indeed a unique solution
with the desired properties. We first consider the {\em homogeneous}
equation obtained by setting $a(r) \equiv 0$ in \eqref{Omexp}. Since
$h(r) \to 1$ as $r \to 0$ and $h(r)$ decays rapidly as $r \to \infty$, 
this linear equation has two particular solutions $\psi_+$, $\psi_-$ 
which satisfy
\begin{equation}\label{psipmdef}
  \psi_-(r) \,\sim\, r^n \quad \hbox{as }r \to 0~, \quad 
  \hbox{and}\quad \psi_+(r) \,\sim\, r^{-n} \quad \hbox{as }r 
  \to \infty~.
\end{equation}
These solutions are of course unique. Moreover, since $n \ge 2$, 
the coefficient $n^2/r^2 - h(r)$ in \eqref{Omexp} is always
positive, because $\sup_{r > 0}r^2 h(r) \cong 2.59\ldots < n^2$. 
It then follows from the Maximum Principle that the functions 
$\psi_+$, $\psi_-$ are strictly positive and satisfy
\[
  \psi_-(r) \,\sim\, \kappa_- r^n \quad \hbox{as }r \to \infty~, \quad 
  \hbox{and}\quad \psi_+(r) \,\sim\, \kappa_+ r^{-n} \quad \hbox{as }r 
  \to 0~,
\]
for some $\kappa_-, \kappa_+ > 0$. In particular, $\psi_+$, $\psi_-$
are linearly independent. In addition, their Wronskian determinant 
satisfies $W = \psi_+ \psi_-' - \psi_- \psi_+' = w_0/r$ for some 
$w_0 > 0$, and it follows that $\kappa_+ = \kappa_- = w_0/(2n)$. 
We now return to the full equation \eqref{Omexp}, and consider the 
particular solution given by the explicit formula
\begin{equation}\label{Omsol}
  \Omega(r) \,=\, \psi_+(r)\int_0^r \frac{y}{w_0}\,\psi_-(y) 
  \frac{a(y)}{n\phi(y)}\d y + \psi_-(r)\int_r^\infty \frac{y}{w_0}
  \,\psi_+(y)\frac{a(y)}{n\phi(y)}\d y~, \quad r > 0~. 
\end{equation}
As is easily verified, we have $\Omega(r) = \cO(r^n)$ as $r \to 0$,  
$\Omega(r) = \cO(r^{-n})$ as $r \to \infty$, and $\Omega$ is the unique 
solution of \eqref{Omexp} with these properties. If $w = -\omega(r)
\cos(n\theta)$, where $\omega$ is defined by \eqref{omexp}, then 
$\Lambda w = z$ by construction, and it is not difficult 
to see that $w \in Y_n \cap Z$. Indeed, it is clear that $w$ is 
smooth away from the origin, and the fact that $e^{|\xi|^2/8}w$ 
decays rapidly at infinity follows immediately from \eqref{omexp}. 
To prove that $w$ is smooth in a neighborhood of zero, we observe 
that any regular solution of \eqref{Omexp} has the form $\Omega(r)
 = r^n \Phi(r^2)$, where $\Phi : [0,\infty) \to \real$ is a smooth
function. Using \eqref{omexp}, we conclude that
\[
  w(\xi) \,=\, -|\xi|^n \cos(n\theta) \Bigl(h(|\xi|)\Phi(|\xi|^2)
  + \frac{A(|\xi|^2)}{n\phi(|\xi|)}\Bigr)
\]
is smooth also near the origin, because $|\xi|^n \cos(n\theta) = 
\Re((\xi_1 + i\xi_2)^n)$ is a homogeneous polynomial in $\xi$. 
The proof of Lemma~\ref{Laminv} is thus complete. \QED

\medskip\noindent{\bf Remark.} Of course, the argument above 
fails if $n = 1$, because the coefficient $1/r^2 - h(r)$ in 
\eqref{Omexp} is no longer positive. In fact, it is easy to 
verify that the functions $\psi_+, \psi_-$ defined by 
\eqref{psipmdef} are linearly dependent in that case. 

\medskip
Equipped with Lemma~\ref{Laminv}, we now go back to the determination
of the vorticity profile $F_i(\xi,t)$ in \eqref{wvfirst}. The 
equation we have to solve is $\Lambda F_i + A_i = 0$, where 
$A_i(\xi,t)$ is given by \eqref{ABC}. Using polar coordinates 
$(r,\theta)$ as before, we can write
\[
  A_i(\cdot,t) \,=\, \frac{d^2}{4\pi}\,r^2g(r)\sum_{j \neq i}
  \frac{\alpha_j}{\alpha_i}\,\frac{1}{|z_{ij}(t)|^2}
  \,\sin(2(\theta-\theta_{ij}(t)))~,
\]
where $\theta_{ij}(t)$ is the argument of the vector $z_{ij}(t) =
z_i(t)-z_j(t)$. This expression shows that $A_i(\cdot,t) \in Y_2 \cap
Z$ for any $t \in [0,T]$, hence by Lemma~\ref{Laminv} there exists a 
unique $F_i^0(\cdot,t) \in Y_2 \cap Z$ such that $\Lambda F_i^0 + 
A_i = 0$. Explicitly,
\begin{equation}\label{Fi0def}
  F_i^0(\cdot,t) \,=\, \frac{d^2}{4\pi}\,\omega(r)\sum_{j \neq i}
  \frac{\alpha_j}{\alpha_i}\,\frac{1}{|z_{ij}(t)|^2}
  \,\cos(2(\theta-\theta_{ij}(t)))~,
\end{equation}
where $\omega(r)$ is given by \eqref{omexp}, \eqref{Omexp} 
with $n = 2$ and $a(r) = r^2 g(r)$. It follows in particular from 
\eqref{omexp}, \eqref{Fi0def} that 
\[
  |F_i^0(\xi,t)| \,\le\, C|\xi|^2(1+|\xi|^2)\,e^{-|\xi|^2/4}~, \qquad
  \xi \in \real^2~, \quad t \in (0,T]~.
\]
If we now return to \eqref{R1def} and choose $F_i(\xi,t) = F_i^0(\xi,t)$, 
it is easy to verify that
\[
  R_i^{(1)}(\xi,t) \,=\, R_i^{(0)}(\xi,t) - \frac{\alpha_i t}{d^2}
  \,A_i(\xi,t) + \tilde R_i^{(1)}(\xi,t) \,=\, 
  \cO\Bigl(\Bigl(\frac{\nu t}{d^2}\Bigr)^{\frac12}\Bigr)~,  
\]
so we succeded in constructing an approximate solution with a smaller
residuum than $R_i^{(0)}$. As is explained in Section~\ref{sec2}, the
profile $F_i^0(\xi,t)$ describes to leading order the deformations of
the vortices due to mutual interaction.

\medskip
Before going further, we state and prove a variant of 
Lemma~\ref{Laminv} which will be useful in the next section.

\begin{lemma}\label{Laminv2}
Assume that $z \in Y_n \cap Z$ for some $n \ge 2$, and let $w \in Y_n 
\cap Z$ be the solution of $\Lambda w = z$ given by Lemma~\ref{Laminv}.
Then for all $\epsilon \neq 0$ the equation
\begin{equation}\label{epsLam}
  \epsilon(1 - \cL)w^\epsilon + \Lambda w^\epsilon \,=\, z
\end{equation}
has a unique solution $w^\epsilon \in Y_n \cap Z$. Moreover, there 
exists $C > 0$ (depending on $z$) such that
\begin{equation}\label{epscomp}
  \|w^\epsilon - w\|_Y \,\le\, \frac{C|\epsilon|}{1+|\epsilon|}~,
  \quad \hbox{for all } \epsilon \neq 0~. 
\end{equation}
\end{lemma}

\medskip\noindent{\bf Proof.} Again, a particular case of 
Lemma~\ref{Laminv2} has been proved in \cite[Proposition~3.4]{GW07}. 
As is well-known (see e.g. \cite{GW02}), the operator $\cL$ defined
by \eqref{cLdef} is self-adjoint in $Y$ and its spectrum is given 
by $\sigma(\cL) = \{-\frac{n}{2}\,|\, n = 0,1,2,\dots\}$. Since 
$\Lambda$ is skew-symmetric and relatively compact with respect to
$\cL$, it follows that the operator $-\epsilon\cL + \Lambda$ is 
maximal accretive for any $\epsilon > 0$. Thus, for any $z \in Y$, 
equation \eqref{epsLam} has a unique solution $w^\epsilon \in Y$, 
which satisfies $\|w^\epsilon\|_Y \le \epsilon^{-1}\|z\|_Y$. A 
similar result holds of course for $\epsilon < 0$. 

We now consider the particular case where $z \in Y_n \cap Z$ for 
some $n \ge 2$. Since $\cL$ commutes with the projection $P_n$,
it is clear that $P_n w^\epsilon = w^\epsilon$, hence $w^\epsilon 
\in Y_n$. In that subspace, Eq.~\eqref{epsLam} reduces to an ordinary
differential equation which can be studied as in the proof of 
Lemma~\ref{Laminv}. In particular, it is straightforward to check 
that $e^{|\xi|^2/8}w^\epsilon$ decays rapidly as $|\xi| \to \infty$, 
so that $w^\epsilon \in Y_n \cap Z$. On the other hand, since
\begin{equation}\label{Lamaux1}
  \epsilon(1-\cL)(w^\epsilon-w) + \Lambda(w^\epsilon-w) \,=\, 
  -\epsilon(1-\cL)w~,
\end{equation}
the fact that $\|[\epsilon(1-\cL) + \Lambda]^{-1}\| \le |\epsilon|^{-1}$ 
implies
\begin{equation}\label{Lamaux2}
  \|w^\epsilon-w\|_Y \,\le\, \|[\epsilon(1-\cL) + \Lambda]^{-1}\| 
  \,\|\epsilon (1-\cL)w\|_Y \,\le\, \|(1-\cL)w\|_Y~,
\end{equation}
hence $\|w^\epsilon-w\|_Y$ is uniformly bounded for all $\epsilon 
\neq 0$. 

Finally, since $w \in Y_n \cap Z$, we have $(1-\cL)w \in Y_n \cap Z$,
hence by Lemma~\ref{Laminv} there exists a unique $\hat w \in Y_n 
\cap Z$ such that $\Lambda \hat w = (1-\cL)w$. Then \eqref{Lamaux1} 
takes the equivalent form
\[
  [\epsilon(1-\cL) + \Lambda](w^\epsilon-w + \epsilon\hat w) \,=\,
  \epsilon^2 (1-\cL)\hat w~,
\]
from which we deduce
\[
  w^\epsilon - w \,=\, -\epsilon\hat w + \epsilon^2 
  [\epsilon(1-\cL) + \Lambda]^{-1} (1-\cL)\hat w~.
\]
Since $\|[\epsilon(1-\cL) + \Lambda]^{-1}\| \le |\epsilon|^{-1}$, 
we conclude that
\begin{equation}\label{Lamaux3}
  \|w^\epsilon-w\|_Y \,\le\, |\epsilon|(\|\hat w\|_Y + \|(1-\cL)
  \hat w\|_Y)~. 
\end{equation}
Combining \eqref{Lamaux2} and \eqref{Lamaux3}, we obtain \eqref{epscomp}. 
This concludes the proof of Lemma~\ref{Laminv2}. \QED

\begin{remark}\label{Zstar}
For later use, we also introduce the space $Z_* \subset Z$ defined by
\[
  Z_* \,=\, \Bigl\{w : \real^2 \to \real\,\Big|\, e^{|\xi|^2/4} w \in 
  \cS_*(\real^2)\Bigr\}~,
\]
where $\cS_*(\real^2)$ denotes the space of all smooth functions $w :
\real^2 \to \real$ such that $w$ and all its derivatives have at 
most a polynomial growth at infinity. As is easily verified,
Lemmas~\ref{Laminv} and \ref{Laminv2} still hold if we replace
everywhere $Z$ by $Z_*$. In particular, if $z \in Y_n \cap Z_*$,
the solution $w^\epsilon$ of \eqref{epsLam} belongs to $Y_n \cap
Z_*$. Thus, for any $\gamma < 1$, there exists $C > 0$ such that
\[
  |w^\epsilon(\xi)| + |\nabla w^\epsilon(\xi)| \,\le\, C\,
  e^{-\gamma|\xi|^2/4}~, \qquad \xi \in \real^2~,
\]
and using for instance \eqref{Lamaux1} one can show that the constant 
$C$ is independent of $\epsilon$. 
\end{remark}

\subsection{Third order approximation}\label{3.3}

We now construct an approximate solution that will be accurate 
enough to prove Theorem~\ref{thm4}. We set
\begin{equation}\label{wvapp}
\begin{array}{rcl}
  w_i^\app(\xi,t) &=& {\DS G(\xi) + \Bigl(\frac{\nu t}{d^2}\Bigr)
  F_i(\xi,t) + \Bigl(\frac{\nu t}{d^2}\Bigr)^{3/2}H_i(\xi,t) + 
  \Bigl(\frac{\nu t}{d^2}\Bigr)^2 K_i(\xi,t)~,} \\ [4mm]
  v_i^\app(\xi,t) &=& {\DS v^G(\xi) + \Bigl(\frac{\nu t}{d^2}\Bigr)
  v^{F_i}(\xi,t) + \Bigl(\frac{\nu t}{d^2}\Bigr)^{3/2}v^{H_i}(\xi,t) + 
  \Bigl(\frac{\nu t}{d^2}\Bigr)^2 v^{K_i}(\xi,t)~,}
\end{array}
\end{equation}
where the vorticity profiles $F_i$, $H_i$, $K_i$ have to be determined, 
and the velocity fields $v^{F_i}$, $v^{H_i}$, $v^{K_i}$ are obtained 
from $F_i$, $H_i$, $K_i$ via the Biot-Savart law \eqref{BS}. As in 
Proposition~\ref{R0expand}, the main expansion parameter in 
\eqref{wvapp} is $(\nu t/d^2)$. The profiles $F_i$, $H_i$, $K_i$ still 
depend on the viscosity $\nu$, but they all have a finite limit as 
$\nu \to 0$. 

Our first task is to compute the residuum $R_i^{(3)}(\xi,t)$ of
the third-order expansion \eqref{wvapp}, as an approximate solution
of \eqref{wieq}. By a direct calculation, we find
\begin{align}\nonumber
  (t\partial_t -\cL)w_i^\app(\xi,t) \,&=\, \Bigl(\frac{\nu t}{d^2}\Bigr)
  \Bigl(t\partial_t F_i + F_i -\cL F_i\Bigr)(\xi,t) \\ \label{res1}
  \,&+\, \Bigl(\frac{\nu t}{d^2}\Bigr)^{3/2}\Bigl(t\partial_t H_i + 
  \frac32 H_i -\cL H_i\Bigr)(\xi,t) \\ \nonumber
  \,&+\, \Bigl(\frac{\nu t}{d^2}\Bigr)^2
  \Bigl(t\partial_t K_i + 2K_i -\cL K_i\Bigr)(\xi,t)~,
\end{align}
and using \eqref{R0def2}, \eqref{Lamdef} we obtain
\begin{align}\nonumber
\sum_{j=1}^N &\frac{\alpha_j}{\nu}\Bigl\{v_j^\app\Bigl(\xi + 
  \frac{z_{ij}(t)}{\sqrt{\nu t}}\,,t\Bigr) - v^G
  \Bigl(\frac{z_{ij}(t)}{\sqrt{\nu t}}\Bigr)\Bigr\}\cdot 
  \nabla w_i^\app(\xi,t) \\ \nonumber
\,&=\, R_i^{(0)}(\xi,t) + \frac{\alpha_i t}{d^2}\Bigl\{\Lambda F_i 
  + \Bigl(\frac{\nu t}{d^2}\Bigr)^{\frac12}\Lambda H_i + 
  \Bigl(\frac{\nu t}{d^2}\Bigr)\Lambda K_i\Bigr\}(\xi,t) \\ 
  \label{res2}
\,&+\,\sum_{j\neq i} \frac{\alpha_j}{\nu}\Bigl\{
  \Bigl(\frac{\nu t}{d^2}\Bigr)v^{F_j} +  \Bigl(\frac{\nu t}{d^2}
  \Bigr)^{\frac32}v^{H_j} + \Bigl(\frac{\nu t}{d^2}\Bigr)^2v^{K_j}\Bigr\}
  \Bigl(\xi + \frac{z_{ij}(t)}{\sqrt{\nu t}}\,,t\Bigr) \cdot
  \nabla G(\xi) \\ \nonumber
\,&+\,\sum_{j\neq i} \frac{\alpha_j}{\nu}\Bigl\{v^G\Bigl(\xi + 
  \frac{z_{ij}(t)}{\sqrt{\nu t}}\Bigr) - v^G
  \Bigl(\frac{z_{ij}(t)}{\sqrt{\nu t}}\Bigr)\Bigr\}\cdot 
  \nabla \Bigl\{\Bigl(\frac{\nu t}{d^2}\Bigr)F_i +  \Bigl(\frac{\nu t}{d^2}
  \Bigr)^{\frac32}H_i + \Bigl(\frac{\nu t}{d^2}\Bigr)^2 K_i\Bigr\}(\xi,t) 
  \\ \nonumber
\,&+\,\sum_{j=1}^N \frac{\alpha_j}{\nu}\Bigl\{
  \Bigl(\frac{\nu t}{d^2}\Bigr)v^{F_j} +  \Bigl(\frac{\nu t}{d^2}
  \Bigr)^{\frac32}v^{H_j} + \Bigl(\frac{\nu t}{d^2}\Bigr)^2v^{K_j}\Bigr\}
  \Bigl(\xi + \frac{z_{ij}(t)}{\sqrt{\nu t}}\,,t\Bigr)\,\cdot \\ \nonumber
\,& \qquad \quad\cdot\, \nabla \Bigl\{\Bigl(\frac{\nu t}{d^2}\Bigr)F_i + 
  \Bigl(\frac{\nu t}{d^2} \Bigr)^{\frac32}H_i + \Bigl(\frac{\nu t}{d^2}
  \Bigr)^2 K_i\Bigr\}(\xi,t)~.
\end{align}
By definition, the residuum $R_i^{(3)}(\xi,t)$ is the sum of all 
terms in \eqref{res1} and \eqref{res2}. 

Next, we separate the lower order terms, which will have to be eliminated 
by an appropriate choice of $F_i, H_i, K_i$, from the higher order 
terms, which will be automatically negligible. We thus decompose
\[
  R_i^{(3)}(\xi,t) \,=\, R_i^\ell(\xi,t) + R_i^h(\xi,t)~, 
\]
where $R_i^\ell$ collects the lower order terms in the residuum, 
namely
\begin{align}\nonumber
  R_i^\ell(\xi,t) \,&=\, 
    R_i^{(0)}(\xi,t) + \frac{\alpha_i t}{d^2}\Bigl\{\Lambda F_i 
  + \Bigl(\frac{\nu t}{d^2}\Bigr)^{\frac12}\Lambda H_i + 
  \Bigl(\frac{\nu t}{d^2}\Bigr)\Lambda K_i\Bigr\}(\xi,t) \\ \label{lowdef}
  \,&+\, \Bigl(\frac{\nu t}{d^2}\Bigr) \Bigl(t\partial_t F_i + F_i 
    -\cL F_i\Bigr)(\xi,t) + \frac{\alpha_i t}{d^2}\Bigl(\frac{\nu t}
    {d^2}\Bigr)v^{F_i}(\xi,t)\cdot\nabla F_i(\xi,t) \\ \nonumber
  \,&+\, \Bigl(\frac{\nu t}{d^2}\Bigr) \sum_{j\neq i} \frac{\alpha_j}
    {\nu}\Bigl\{v^G\Bigl(\xi + \frac{z_{ij}(t)}{\sqrt{\nu t}}\Bigr) 
    - v^G\Bigl(\frac{z_{ij}(t)}{\sqrt{\nu t}}\Bigr)\Bigr\}\cdot 
    \nabla F_i(\xi,t)~.
\end{align}

Our goal is to choose the vorticity profiles $F_i$, $H_i$, $K_i$ so as
to minimize the quantity $R_i^\ell(\xi,t)$ in the vanishing viscosity
limit. The contributions of the naive residuum \eqref{R0exp} are
easily eliminated by successive applications of Lemma~\ref{Laminv}: we
first take $F_i(\xi,t)$ so that $\Lambda F_i + A_i = 0$, then
$H_i(\xi,t)$ so that $\Lambda H_i + B_i = 0$, and so on. In this way,
we obtain a residuum $R_i^\ell(\xi,t)$ of size $\cO(\nu t/d^2)$, but
unfortunately this is not sufficient to prove Theorems~\ref{thm3} and
\ref{thm4}. To obtain a more precise estimate, we also need to
eliminate all terms in the last two lines of \eqref{lowdef}. This
requires a more careful choice of $F_i$, $H_i$, $K_i$, which will be
done in three steps:

\medskip\noindent{\bf 1)} First, we take
\begin{equation}\label{Fidef}
  F_i(\xi,t) \,=\, \bar F_i(\xi,t) + F_i^\nu(\xi,t)~, \qquad 
  \xi \in \real^2~, \quad t \in (0,T]~,
\end{equation}
where $\bar F_i(\cdot,t) \in Y_0 \cap Z$ and $F_i^\nu(\cdot,t)
\in Y_2 \cap Z$. More precisely: 

\smallskip\noindent{\sl a)} The profile $F_i^\nu(\xi,t)$ is the 
unique solution of the elliptic equation
\begin{equation}\label{Finudef}
  \frac{\nu}{\alpha_i}(F_i^\nu - \cL F_i^\nu) + \Lambda F_i^\nu 
  + A_i \,=\, 0~,
\end{equation}
as given by Lemma~\ref{Laminv2}. Since $A_i(\cdot,t) \in Y_2 \cap Z$, 
we have $F_i^\nu(\cdot,t) \in Y_2 \cap Z$ for all $t \in (0,T]$. 
Moreover, by \eqref{epscomp}, 
\[
  \|F_i^\nu(\cdot,t) - F_i^0(\cdot,t)\|_Y \,\le\, C\,\frac{\nu}
  {|\alpha_i|+\nu}~, \qquad t \in (0,T]~,
\]
where $F_i^0(\xi,t)$ is given by \eqref{Fi0def}. 

\smallskip\noindent{\sl b)} The profile $\bar F_i(\xi,t)$ is the 
unique solution of the linear parabolic equation
\begin{equation}\label{Fibardef}
  t\partial_t \bar F_i + \bar F_i -\cL \bar F_i + \frac{\alpha_i t}{d^2}
  \Bigl\{P_0(V^{F_i^\nu}\cdot\nabla F_i^\nu) + P_0(D_i \cdot \nabla F_i^\nu)
  \Bigr\} \,=\, 0~,
\end{equation}
with initial data $\bar F_i(\cdot,0) = 0$. Here $P_0$ is the orthogonal
projection (in $Y$) onto the radially symmetric functions, and $D_i(\xi,t)$
is the divergence-free vector field given by
\begin{equation}\label{Didef}
  D_i(\xi,t) \,=\, \frac{1}{2\pi}\sum_{j\neq i}\frac{\alpha_j}{\alpha_i}\,
  \frac{d^2}{|z_{ij}(t)|^4}\Bigl(\xi^\perp |z_{ij}(t)|^2 - 2
  (\xi\cdot z_{ij}(t)) z_{ij}(t)^\perp\Bigr)~.
\end{equation}
It is clear that $\bar F_i(\xi,t)$ is well-defined. In fact, if 
$S(\tau) = \exp(\tau\cL)$ denotes the $C_0$-semigroup in $Y$ generated
by $\cL$, we have the explicit formula
\begin{equation}\label{Fibarexp}
  \bar F_i(\cdot,t) \,=\, -\frac{\alpha_i}{d^2}\int_0^t 
  S\Bigl(\log\frac{t}{s}\Bigr)\frac{s}{t} \,P_0 Q_i(\cdot,s)\d s~,
\end{equation}
where $Q_i = V^{F_i^\nu}\cdot\nabla F_i^\nu + D_i \cdot \nabla F_i^\nu$. 

\medskip\noindent{\bf 2)}
Next, we determine the profile $H_i(\xi,t)$. From the proof of 
Proposition~\ref{R0expand}, it is clear that $B_i(\cdot,t) \in Y_3 
\cap Z$. Thus, by Lemma~\ref{Laminv}, there exists a unique solution 
$H_i(\cdot,t) \in Y_3 \cap Z$ to the equation
\begin{equation}\label{Hidef}
  \Lambda H_i(\cdot,t) + B_i(\cdot,t) \,=\, 0~, \qquad t \in (0,T]~. 
\end{equation}

\medskip\noindent{\bf 3)} Finally, we set
\begin{equation}\label{Kidef}
  K_i(\xi,t) \,=\, K_{i2}(\xi,t) + K_{i4}(\xi,t)~, \qquad 
  \xi \in \real^2~, \quad t \in (0,T]~,
\end{equation}
where $K_{i2}(\cdot,t) \in Y_2 \cap Z$ and $K_{i4}(\cdot,t) 
\in Y_4 \cap Z$ are chosen as follows:

\smallskip\noindent{\sl a)} The profile $K_{i2}(\xi,t)$ is the 
unique solution, given by Lemma~\ref{Laminv}, of the equation
\begin{equation}\label{Ki2def}
  \Lambda K_{i2} + P_2 (V^{F_i}\cdot\nabla F_i) + P_2(D_i \cdot 
  \nabla F_i) + \frac{d^2}{\alpha_i}\partial_t F_i^\nu \,=\, 0~, 
\end{equation}
where $P_2$ is the orthogonal projection in $Y$ onto $Y_2$. 

\smallskip\noindent{\sl b)} The profile $K_{i4}(\xi,t)$ is the 
unique solution, given by Lemma~\ref{Laminv}, of the equation
\begin{equation}\label{Ki4def}
  \Lambda K_{i4} + P_4 (V^{F_i}\cdot\nabla F_i) + P_4(D_i \cdot 
  \nabla F_i) + C_i \,=\, 0~, 
\end{equation}
where $P_4$ is the orthogonal projection onto $Y_4$, and $C_i$ is
as in \eqref{ABC}. 

\medskip The main result of this section is:

\begin{proposition}\label{residuum}
Fix $\frac12 < \gamma < 1$. There exists $C > 0$ such that, 
with the above choices of the vorticity profiles $F_i$, $H_i$, 
$K_i$, the residuum of the approximate solution \eqref{wvapp}
satisfies
\begin{equation}\label{Ri3}
  |R_i^{(3)}(\xi,t)| \,\le\,  C \,\Bigl(\frac{\nu t}{d^2}\Bigr)^{3/2}
  \,e^{-\gamma|\xi|^2/4}~, 
\end{equation}
for all $\xi \in \real^2$, all $t \in (0,T]$, and all $i \in 
\{1,\cdots,N\}$. 
\end{proposition}

\noindent{\bf Proof.} The proof is a long sequence of rather 
straightforward verifications, some of which will be left to the 
reader. We first summarize the informations we have on the vorticity
profiles $F_i$, $H_i$, $K_i$, and on the associated velocity fields
$v^{F_i}$, $v^{H_i}$, $v^{K_i}$. Given any $\gamma < 1$, we claim 
that there exists $C > 0$ such that
\begin{equation}\label{Fibdd}
  |F_i(\xi,t)| + |\nabla F_i(\xi,t)| \,\le\, C\,e^{-\gamma|\xi|^2/4}~,
  \qquad |v^{F_i}(\xi,t)| \,\le\, \frac{C}{(1+|\xi|^2)^{3/2}}~, 
\end{equation}
for all $\xi \in \real^2$, all $t \in (0,T]$, and all $i \in
\{1,\cdots,N\}$. Indeed, we recall that $F_i = \bar F_i + F_i^\nu$,
where $F_i^\nu$ is defined by \eqref{Finudef} and $\bar F_i$ by
\eqref{Fibardef}. From \eqref{Finudef} and Remark~\ref{Zstar}, we know
that $F_i^\nu(\cdot,t) \in Y_2 \cap Z_*$. Thus $F_i^\nu$ and $\nabla
F_i^\nu$ satisfy the Gaussian bound in \eqref{Fibdd}, and it follows
from \cite[Proposition~B.1]{GW02} that $|v^{F_i}(\xi,t)| \le
C(1+|\xi|^2)^{-3/2}$, see also \eqref{vnbdd}. On the other hand, using
\eqref{Fibarexp} and the explicit expression of the integral kernel of
the semigroup $S(\tau) = \exp(\tau\cL)$ \cite[Appendix~A]{GW02}, 
it is straightforward to verify that $|\bar F_i(\xi,t)| \le C 
e^{-\gamma|\xi|^2/4}$. Since $\nabla S(\tau) \,=\, e^{\tau/2}S(\tau)
\nabla$, we also have
\[
  \nabla \bar F_i(\cdot,t) \,=\, -\frac{\alpha_i}{d^2}\int_0^t 
  S\Bigl(\log\frac{t}{s}\Bigr)\Bigl(\frac{s}{t}\Bigr)^{\frac12} 
  \,\nabla P_0 Q_i(\cdot,s)\d s~,
\]
from which we deduce that $|\nabla \bar F_i(\xi,t)| \le C 
e^{-\gamma|\xi|^2/4}$. Finally, since $\bar F_i(\cdot,t)$ is 
a radially symmetric function with zero average, we have a simple
formula for the associated velocity field
\[
  v^{\bar F_i}(\xi,t) \,=\, \frac{1}{2\pi}\,\frac{\xi^\perp}{|\xi|^2}
  \int_{|\xi'| \ge |\xi|} \bar F_i(\xi',t)\d\xi'~,
\]
which implies that $v^{\bar F_i}(\xi,t)$ has also a Gaussian decay
as $|\xi| \to \infty$. This proves \eqref{Fibdd}. 

The corresponding estimates for $H_i(\xi,t)$ and $K_i(\xi,t)$ are
easier to establish. Since $B_i(\cdot,t) \in Y_3 \cap Z_*$, 
it follows from \eqref{Hidef} and Remark~\ref{Zstar} that 
$H_i(\cdot,t) \in Y_3 \cap Z_*$ for all $t \in (0,T]$. Using the
same arguments as before, we obtain
\begin{equation}\label{Hibdd}
  |H_i(\xi,t)| + |\nabla H_i(\xi,t)| \,\le\, C\,e^{-\gamma|\xi|^2/4}~,
  \qquad |v^{H_i}(\xi,t)| \,\le\, \frac{C}{(1+|\xi|^2)^2}~. 
\end{equation}
Similarly, it follows from \eqref{Ki2def}, \eqref{Ki4def} that 
$K_{i2}(\cdot,t) \in Y_2 \cap Z_*$ and $K_{i4}(\cdot,t) \in Y_4 
\cap Z_*$. We conclude that $K_i = K_{i2} + K_{i4}$ satisfies 
\begin{equation}\label{Kibdd}
  |K_i(\xi,t)| + |\nabla K_i(\xi,t)| \,\le\, C\,e^{-\gamma|\xi|^2/4}~,
  \qquad |v^{K_i}(\xi,t)| \,\le\, \frac{C}{(1+|\xi|^2)^{3/2}}~. 
\end{equation}

Next, we make the following observations, which were implicitely 
used in the definitions of the profiles $F_i$ and $K_i$. Since 
$F_i = \bar F_i + F_i^\nu \in Y_0 + Y_2$, it is easy to verify that
$v^{F_i}\cdot\nabla F_i \in Y_0 + Y_2 + Y_4$. Similarly, using 
the definition \eqref{Didef} of the vector field $D_i$, we find 
that $D_i \cdot\nabla F_i \in Y_0 + Y_2 + Y_4$. Thus we have 
the identities
\begin{equation}\label{Idd1}
  v^{F_i}\cdot\nabla F_i \,=\, (P_0 + P_2 + P_4)(v^{F_i}\cdot\nabla F_i)~,
  \qquad D_i\cdot\nabla F_i \,=\, (P_0 + P_2 + P_4)(D_i\cdot\nabla F_i)~.
\end{equation}
Moreover, it is straightforward to check that
\begin{equation}\label{Idd2}
  P_0(v^{F_i}\cdot\nabla F_i) \,=\, P_0(v^{F_i^\nu}\cdot\nabla F_i^\nu)~, 
  \qquad P_0(D_i\cdot\nabla F_i) \,=\, P_0(D_i \cdot\nabla F_i^\nu)~.
\end{equation}
Remark that both expressions in \eqref{Idd2} appear in the definition
\eqref{Fibardef} of $\bar F_i$. 

Now, we replace the definitions \eqref{Fidef}--\eqref{Fibardef} and
\eqref{Hidef}--\eqref{Ki4def} into the expression \eqref{lowdef}
of the residuum $R_i^\ell(\xi,t)$. Using in addition \eqref{R0exp}, 
\eqref{Idd1}, \eqref{Idd2}, we obtain the simple formula
\begin{equation}\label{lowexp}
  R_i^\ell(\xi,t) \,=\, \frac{\alpha_i t}{d^2}\,\tilde R_i^{(0)}(\xi,t)
  + \Bigl(\frac{\nu t}{d^2}\Bigr) \,\Delta_i(\xi,t)\cdot \nabla
  F_i(\xi,t)~,
\end{equation}
where
\begin{equation}\label{Deltaidef}
  \Delta_i(\xi,t) \,=\, \sum_{j\neq i} \frac{\alpha_j}
  {\nu}\Bigl(v^G\Bigl(\xi + \frac{z_{ij}(t)}{\sqrt{\nu t}}\Bigr) 
   - v^G\Bigl(\frac{z_{ij}(t)}{\sqrt{\nu t}}\Bigr)\Bigr) - 
    \frac{\alpha_i t}{d^2}\,D_i(\xi,t)~.
\end{equation}
Our goal is to obtain an estimate of the form \eqref{Ri3} for 
$R_i^\ell(\xi,t)$. Since $\tilde R_i^{(0)}(\xi,t)$ satisfies 
\eqref{R0rem}, it is sufficient to bound the second term in 
the right-hand side of \eqref{lowexp}. As in the proof of 
Proposition~\ref{R0expand}, we can assume that $|\xi| \le 
d/(2\sqrt{\nu t})$, because in the converse case the quantity 
$|\nabla F_i(\xi,t)|$ is extremely small due to \eqref{Fibdd}. 
Using the decomposition \eqref{V1V2} and the definition \eqref{Didef}
of the vector field $D_i$, we find
\begin{align}\nonumber
  \Delta_i(\xi,t) \,&=\, \frac{1}{2\pi}\sum_{j\neq i} \frac{\alpha_j}{\nu}
  \Bigl(V_1\Bigl(\xi,\frac{z_{ij}(t)}{\sqrt{\nu t}}\Bigr) + 
  V_2\Bigl(\xi,\frac{z_{ij}(t)}{\sqrt{\nu t}}\Bigr)\Bigr) - 
  \frac{\alpha_i t}{d^2}\,D_i(\xi,t) \\ \label{Deltaiexp}
  \,&=\, \frac{1}{2\pi}\sum_{j\neq i} \frac{\alpha_j}{\nu}
  W\Bigl(\xi,\frac{z_{ij}(t)}{\sqrt{\nu t}}\Bigr) +
  \frac{1}{2\pi}\sum_{j\neq i} \frac{\alpha_j}{\nu} V_2\Bigl(\xi,
  \frac{z_{ij}(t)}{\sqrt{\nu t}}\Bigr)~,
\end{align}
where
\begin{equation}\label{Wexp}
  W(\xi,\eta) \,=\, \frac{(\xi+\eta)^\perp}{|\xi+\eta|^2} 
  - \frac{(\xi+\eta)^\perp}{|\eta|^2} + 2\frac{(\xi\cdot\eta)\eta^\perp}
  {|\eta|^4} \,=\, \cO\Bigl(\frac{|\xi|^2}{|\eta|^3}\Bigr)~, 
  \quad \hbox{as } |\eta| \to \infty~.    
\end{equation}
In view of \eqref{V2exp}, the contributions of $V_2$ are negligible, 
and using \eqref{Fibdd} we easily obtain
\[
  |\Delta_i(\xi,t) \cdot \nabla F_i(\xi,t)| \,\le\, 
  C \frac{|\alpha_i| t}{d^2}\Bigl(\frac{\nu t}{d^2}\Bigr)^{1/2}
  e^{-\gamma|\xi|^2/4}~, \qquad \xi \in \real^2~, 
\]
which is the desired estimate. 

To complete the proof of Proposition~\ref{residuum}, it remains 
to verify that 
\[
  |R_i^{(3)}(\xi,t) - R_i^\ell(\xi,t)| \,\le\, C
  \,\Bigl(\frac{\nu t}{d^2}\Bigr)^{3/2}\,e^{-\gamma|\xi|^2/4}~.
\]
This follows immediately from \eqref{res2}, \eqref{lowdef} if one
uses the bounds \eqref{Fibdd}, \eqref{Hibdd}, \eqref{Kibdd} on
the vorticity profiles $F_i$, $H_i$, $K_i$ and the associated 
velocities. In particular, it is straightforward to check that 
each term in the difference $R_i^{(3)}- R_i^\ell$ is of the order
of $(\nu t/d^2)^n$ for some $n \ge 3/2$, either due to an explicit
prefactor or as a consequence of the polynomial decay of the 
velocity fields $v^G$,  $v^{F_i}$, $v^{H_i}$, or $v^{K_i}$ as 
$|\xi| \to \infty$. This concludes the proof. \QED

\medskip\noindent{\bf Remark.} Instead of \eqref{Hidef}, one can
define the profile $H_i(\xi,t)$ as the (unique) solution of the 
elliptic equation
\[
  \frac{\nu}{\alpha_i}\Bigl(\frac32 H_i - \cL H_i\Bigr) + 
  \Lambda H_i + B_i \,=\, 0~, 
\]
which is the analog of \eqref{Finudef}. In the same spirit, one 
can replace \eqref{Ki2def} by 
\[
  \frac{\nu}{\alpha_i}(2K_{i2} - \cL K_{i2}) + \Lambda K_{i2} 
  + P_2 (V^{F_i}\cdot\nabla F_i) + P_2(D_i \cdot 
  \nabla F_i) + \frac{d^2}{\alpha_i}\partial_t F_i^\nu \,=\, 0~, 
\]
and proceed similarly with \eqref{Ki4def}. After these modifications,
it is easy to verify that the residuum satisfies the improved bound
\[
  |R_i^{(3)}(\xi,t)| \,\le\,  C \,\frac{|\alpha| t}{d^2}
  \,\Bigl(\frac{\nu t}{d^2}\Bigr)^{3/2}\,e^{-\gamma|\xi|^2/4}~, 
\]
which is sharper than \eqref{Ri3} for small times. This refinement 
is not needed in the proof of Theorem~\ref{thm4}, but it indicates 
the correct way to proceed if one wants to construct even more 
precise approximations of the $N$-vortex solution. 


\section{Control of the remainder}\label{sec4}

In the previous section, we constructed an approximate solution
$w_i^\app(\xi,t)$, $v_i^\app(\xi,t)$ of equation \eqref{wieq}, with a
very small residuum. We now consider the exact solution $w_i(\xi,t)$ 
of \eqref{wieq} given by \eqref{wvdef} and Lemma~\ref{lemdec}, and 
we try control the difference $w_i(\xi,t) - w_i^\app(\xi,t)$ for 
$\xi \in \real^2$ and $t \in (0,T]$, in the vanishing viscosity limit. 
To this end, it is convenient to write
\begin{equation}\label{tildewidef}
  w_i(\xi,t) \,=\, w_i^\app(\xi,t) + \Bigl(\frac{\nu t}{d^2}\Bigr)
  \,{\tilde w}_i(\xi,t)~, \quad
  v_i(\xi,t) \,=\, v_i^\app(\xi,t) + \Bigl(\frac{\nu t}{d^2}\Bigr)
  \,{\tilde v}_i(\xi,t)~,
\end{equation}
and to study the evolution system satisfied by the remainder 
$\tilde w_i(\xi,t)$, $\tilde v_i(\xi,t)$. 

Replacing \eqref{tildewidef} into \eqref{wieq}, and using the 
definition \eqref{PW3}, we find
\begin{align}\label{T0}
 t\partial_t \tilde w_i(\xi,t) &- (\cL \tilde w_i)(\xi,t) 
  + \tilde w_i(\xi,t) \\[2mm] 
   \label{T1}
 & + \frac{\alpha_i}{\nu} \Bigl(v_i^\app(\xi,t)\cdot \nabla
   \tilde w_i(\xi,t) + \tilde v_i(\xi,t)\cdot \nabla w_i^\app(\xi,t)
   \Bigr) \\[1mm] \label{T2}
 & + \sum_{j\neq i} \frac{\alpha_j}{\nu}\left\{v_j^\app
   \Bigl(\xi + \frac{z_{ij}(t)}{\sqrt{\nu t}}\,,\,t\Bigr)
   - v^G\Bigl(\frac{z_{ij}(t)}{\sqrt{\nu t}}\Bigr)\right\}\cdot 
   \nabla \tilde w_i(\xi,t) \\ \label{T3}
 & + \sum_{j\neq i} \frac{\alpha_j}{\nu}\,\tilde v_j
   \Bigl(\xi + \frac{z_{ij}(t)}{\sqrt{\nu t}}\,,\,t\Bigr)
   \cdot \nabla w_i^\app(\xi,t) \\ \label{T4}
 & + \sum_{j=1}^N \frac{\alpha_j t}{d^2}\,\tilde v_j
   \Bigl(\xi + \frac{z_{ij}(t)}{\sqrt{\nu t}}\,,\,t\Bigr)
   \cdot \nabla \tilde w_i(\xi,t) + \tilde R_i(\xi,t) \,=\, 0~,  
\end{align}
where $\tilde R_i(\xi,t) = (\nu t/d^2)^{-1}R_i^{(3)}(\xi,t)$. From 
Proposition~\ref{residuum}, we know that
\begin{equation}\label{tildeRi}
  |\tilde R_i(\xi,t)| \,\le\,  C\,\Bigl(\frac{\nu t}{d^2}\Bigr)^{1/2}
  \,e^{-\gamma|\xi|^2/4}~, 
\end{equation}
for all $\xi \in \real^2$, all $t \in (0,T]$, and all 
$i \in \{1,\dots,N\}$. Also, since $\inttwo w_i(\xi,t)\d\xi = 1$
by \eqref{wvdef} and Lemma~\ref{lemdec}, it is clear that 
$\inttwo \tilde w_i(\xi,t)\d\xi = 0$ for all $t \in (0,T]$. 

To prove Theorem~\ref{thm4}, our strategy is to consider the (unique)
solution $\tilde w_i(\xi,t)$ of system~\eqref{T0}--\eqref{T4} with
zero initial data, and to control it on the time interval $(0,T]$
using an energy functional of the form
\begin{equation}\label{Edef}
  E(t) \,=\, \frac12 \sum_{i=1}^N \inttwo p_i(\xi,t) 
  |\tilde w_i(\xi,t)|^2 \d\xi~, \qquad t \in (0,T]~,
\end{equation}
where the weight functions $p_i(\xi,t)$, $i \in \{1,\dots,N\}$,   
will be carefully constructed below. In particular, we shall require 
that $p_i(\xi,t) \ge C e^{\beta |\xi|/4}$ for some $\beta > 0$. 
Using \eqref{T0}--\eqref{tildeRi}, we shall derive a differential 
inequality for $E(t)$ which will imply that $E(t) = \cO(\nu t/d^2)$ 
as $\nu \to 0$. This will show that
\begin{equation}\label{Econc}
  \sum_{i=1}^N \inttwo e^{\beta |\xi|/4} |\tilde w_i(\xi,t)|^2 
  \d\xi \,\le\, C \Bigl(\frac{\nu t}{d^2}\Bigr)~, \qquad t \in (0,T]~,
\end{equation}
for some $C > 0$, if $\nu$ is sufficiently small. In view of 
\eqref{tildewidef}, this estimate is equivalent to \eqref{strongconv2}, 
which is the desired result. 

\subsection{Construction of the weight functions}\label{4.1}

Since the construction of the weights $p_i(\xi,t)$ is rather delicate,
we first explain the main ideas in a heuristic way. Ideally, we would
like to use for each $i \in \{1,\dots,N\}$ the time-independent weight
$p(\xi) = e^{|\xi|^2/4}$, in order to control the remainder $\tilde
w_i(\cdot,t)$ in the function space $Y$ defined in \eqref{Ydef}.  This
is a natural choice for at least two reasons. First, the linear
operator $\cL$ defined in \eqref{cLdef} is self-adjoint in $Y$, and a
straightforward calculation (which will be reproduced below) shows
that
\begin{equation}\label{NegL}
  \inttwo e^{|\xi|^2/4} \tilde w_i(\cL \tilde w_i - \tilde w_i )\d\xi
  \,\le\, -\inttwo e^{|\xi|^2/4} \Bigl(\frac14 |\nabla \tilde 
  w_i|^2 + \frac{|\xi|^2}{24}|\tilde w_i|^2 + \frac12 |\tilde w_i|^2
  \Bigr)\d\xi~.
\end{equation}
Next, using \eqref{wvapp}, we observe that the self-interaction 
terms \eqref{T1} have the form 
\[
  \frac{\alpha_i}{\nu} \Bigl(v_i^\app(\xi,t)\cdot \nabla
  \tilde w_i(\xi,t) + \tilde v_i(\xi,t)\cdot \nabla w_i^\app(\xi,t)
  \Bigr) \,=\, \frac{\alpha_i}{\nu}\,\Lambda \tilde w_i(\xi,t)
  \,+\, \hbox{regular terms}~,
\]
where $\Lambda$ is the linear operator defined in \eqref{Lamdef}.
Here and below, we call ``regular'' all terms which have a finite
limit as $\nu \to 0$. Since $\Lambda$ is skew-symmetric in the space
$Y$ by Proposition~\ref{Lamprop}, we see that the singular term
$(\alpha_i/\nu)\Lambda \tilde w_i$ will not contribute at all to the
variation of the energy if we use the Gaussian weight $e^{|\xi|^2/4}$.

Unfortunately, this naive choice is not appropriate to treat the
advection terms \eqref{T2}, which describe how the perturbation 
$\tilde w_i$ of the $i^{\rm th}$ vortex is transported by the velocity 
field of the other vortices. Indeed, using again \eqref{wvapp}, we 
can write \eqref{T2} as
\begin{equation}\label{T2bis}
  \sum_{j\neq i} \frac{\alpha_j}{\nu}\left\{v^G\Bigl(\xi + 
  \frac{z_{ij}(t)}{\sqrt{\nu t}}\Bigr) - v^G\Bigl(\frac{z_{ij}(t)}
  {\sqrt{\nu t}}\Bigr)\right\}\cdot \nabla \tilde w_i(\xi,t) 
  \,+\, \hbox{regular terms}~.
\end{equation}
Since $v^G(\xi)$ is given by \eqref{Gdef}, and since $|z_{ij}(t)| \ge
d > 0$ when $i \neq j$, it is easy to verify (as in the proof 
of Proposition~\ref{R0expand}) that the first term in \eqref{T2bis} 
has a finite limit as $\nu \to 0$, provided $|\xi| \ll d/\sqrt{\nu t}$. 
Although this is not obvious a priori, we shall see below that the
same term is also harmless if $|\xi| \gg D/\sqrt{\nu t}$, where
\begin{equation}\label{dmax}
  D \,=\, \max_{t \in [0,T]}\,\max_{i\neq j}\,|z_i(t) - z_j(t)|
  \,<\, \infty~. 
\end{equation}
In the intermediate region, however, the first term in \eqref{T2bis}
can be of size $\cO(|\alpha|/\nu)$, and there is no hope to obtain a
better bound. But we should keep in mind that the whole term
\eqref{T2} describes the advection of the perturbation $\tilde w_i$ by
a divergence-free velocity field, and therefore does not contribute to
the variation of the energy in the regions where the weight function
is constant. The idea is thus to modify the Gaussian weight to obtain
a large plateau in the intermediate region where the advection term
\eqref{T2} is singular.

Our second try is therefore a time-dependent weight of the form
\begin{equation}\label{psecond}
  p(\xi,t) \,=\, \left\{\begin{array}{lcl}
  e^{|\xi|^2/4} & \hbox{if} & |\xi| \le \rho(t)~, \\
  e^{\rho(t)^2/4} & \hbox{if} & \rho(t) \le |\xi| \le K\rho(t)~, \\
  e^{|\xi|^2/(4K^2)} & \hbox{if} & |\xi| \ge K\rho(t)~,
  \end{array}\right.  
\end{equation}
where
\[
  \rho(t) \,=\, \frac{d}{2\sqrt{\nu t}}~, \quad \hbox{and} \quad 
  K \,=\, \frac{4D}{d}~.
\]
By construction, the function $\xi \mapsto p(\xi,t)$ coincides with
$e^{|\xi|^2/4}$ in a large disk near the origin, is identically
constant in the intermediate region where the advection terms
\eqref{T2} are dangerous, and becomes Gaussian again when $|\xi|$ is
very large. Moreover, this weight is continuous, radially symmetric,
and satisfies $e^{|\xi|^2/(4K^2)} \le p(\xi,t) \le e^{|\xi|^2/4}$ for
all $\xi \in \real^2$ and all $t \in (0,T]$. Of course, the operator
$\cL$ is no longer self-adjoint in the function space defined by the
modified weight $p(\xi,t)$, but an estimate of the form \eqref{NegL}
can nevertheless be established by a direct calculation. Similarly,
the operator $\Lambda$ is no longer skew-symmetric, but we shall see
below that it remains approximately skew-symmetric with the modified
weight, and this will be sufficient to treat the self-interaction
terms \eqref{T1}.

We now consider the contributions of the advection terms \eqref{T2}
to the variation of the energy \eqref{Edef}, when $p_i(\xi,t) = 
p(\xi,t)$ for all $i \in \{1,\dots,N\}$. Integrating by parts, 
we obtain
\begin{align}\nonumber
  &\inttwo p(\xi,t)\,\tilde w_i(\xi,t) \sum_{j\neq i} 
  \frac{\alpha_j}{\nu} \left\{v_j^\app \Bigl(\xi + \frac{z_{ij}(t)}
  {\sqrt{\nu t}}\,,\,t\Bigr) - v^G\Bigl(\frac{z_{ij}(t)}{\sqrt{\nu t}}
  \Bigr)\right\}\cdot \nabla \tilde w_i(\xi,t) \d\xi \\\label{T2ter}
  & \,=\, -\frac12 \inttwo \,|\tilde w_i(\xi,t)|^2 \sum_{j\neq i}
  \frac{\alpha_j}{\nu} \left\{v_j^\app \Bigl(\xi + \frac{z_{ij}(t)}
  {\sqrt{\nu t}}\,,\,t\Bigr)- v^G\Bigl(\frac{z_{ij}(t)}{\sqrt{\nu t}}
  \Bigr)\right\}\cdot \nabla p(\xi,t)\d\xi~.
\end{align}
If $\rho(t) \le |\xi| \le K\rho(t)$, then $\nabla p(\xi,t) \equiv 0$
by construction, hence it remains to consider the regions where $|\xi|
\le \rho(t)$ or $|\xi| \ge K\rho(t)$. Unfortunately, although the
integrand in \eqref{T2ter} is regular when $|\xi| \le \rho(t)$, the
contributions from that region are still difficult to control if $t$
is large. To see this, we first replace $v_j^\app$ by $v^G$ as in
\eqref{T2bis}, because the difference $v_j^\app - v^G$ is negligible
at this level of analysis.  Assuming that $|\xi| \le \rho(t)$, we find
as in \eqref{Deltaidef}
\[
  \sum_{j\neq i} \frac{\alpha_j}{\nu}\left\{v^G \Bigl(\xi + 
  \frac{z_{ij}(t)}{\sqrt{\nu t}}\Bigr)- v^G\Bigl(\frac{z_{ij}(t)}
  {\sqrt{\nu t}}\Bigr)\right\} \,=\, \frac{\alpha_it}{d^2}
  \left(D_i(\xi,t) + \cO\Bigl(|\xi|^2\frac{\sqrt{\nu t}}{d}\Bigr)
  \right)~,
\]
where $D_i(\xi,t)$ is defined in \eqref{Didef}. As $\nabla p(\xi,t) = 
\frac{\xi}{2}p(\xi,t)$ for $|\xi| \le \rho(t)$, we conclude that 
the main contribution to \eqref{T2ter} has the form
\begin{equation}\label{T2quater}
  \frac{1}{4\pi} \frac{\alpha_i t}{d^2}
  \int_{|\xi| \le \rho(t)} \,|\tilde w_i(\xi,t)|^2 
  \sum_{j\neq i} \frac{\alpha_j}{\alpha_i}\,\frac{d^2 (\xi\cdot z_{ij}(t))
  (\xi\cdot z_{ij}(t)^\perp)}{|z_{ij}(t)|^4}\,p(\xi,t)\d\xi~,
\end{equation}
and is therefore bounded by 
\begin{equation}\label{T2quinquies}
  C\,\frac{|\alpha_i|t}{d^2} \int_{|\xi| \le \rho(t)} 
  |\xi|^2\,e^{|\xi|^2/4} \,|\tilde w_i(\xi,t)|^2 \d\xi~.
\end{equation}
If $t$ is small with respect to the turnover time $T_0$ defined in
\eqref{T0def}, then $|\alpha_i| t \ll d^2$ and the quantity
\eqref{T2quinquies} can be controlled by the negative terms
originating from the diffusion operator $\cL$, see \eqref{NegL}. In
that case, one can show that the expression \eqref{T2ter} is harmless
also in the outer region where $|\xi| \ge K\rho(t)$, and it is not
difficult to verify that the linear terms \eqref{T3} and the nonlinear
terms \eqref{T4} can be controlled in a similar way. Thus, if $T \ll
T_0$, it is possible to carry out the whole proof of
Theorem~\ref{thm4} using the energy functional \eqref{Edef} with
$p_i(\xi,t) = p(\xi,t)$ for all $i \in \{1,\dots,N\}$.

The main difficulty is of course to get rid of the condition $T \ll
T_0$, which is obviously too restrictive. We follow here the same
strategy as in the construction of the approximate solution
$w_i^\app(\xi,t)$ in Sections~\ref{3.2} and \ref{3.3}. So far, we have
used radially symmetric weights to eliminate, or at least to minimize,
the influence of the singular self-interaction terms \eqref{T1}. The
idea is now to add small, nonsymmetric corrections of size $\cO(\nu
t/d^2)$ which, by interacting with the singular expression \eqref{T1},
will produce counter-terms of size $\cO(1)$ that will exactly
compensate for \eqref{T2quater}. Unfortunately, the construction of
these corrections is quite technical, and requires a non-trivial
modification of the underlying radially symmetric weight.

We now give the precise definition of the weights $p_i(\xi,t)$
that will be used in the definition \eqref{Edef} of the energy.
We start from the radially symmetric weight
\begin{equation}\label{p0def}
  p_0(\xi,t) \,=\, \left\{\begin{array}{lcl}
  e^{|\xi|^2/4}\Bigl(1 - \psi(|\xi|^2 - a(t)^2)\Bigr) & \hbox{if} 
  & 0 \le |\xi|^2 \le a(t)^2{+}1~, \\[1mm]
  e^{a(t)^2/4} & \hbox{if} & a(t)^2{+}1 \le |\xi|^2 \le b(t)^2~, \\[1.5mm]
  e^{\beta|\xi|/4} & \hbox{if} & |\xi|^2 \ge b(t)^2~,
  \end{array}\right.  
\end{equation}
where
\begin{equation}\label{abdef}
  a(t) \,=\, a_0 \Bigl(\frac{d^2}{\nu t}\Bigr)^{1/4}~, \quad 
  b(t) \,=\, b_0 \Bigl(\frac{d^2}{\nu t}\Bigr)^{1/2}~, \quad 
 \hbox{and} \quad \beta \,=\, \frac{a_0^2}{b_0} \,\ll\, 1~.
\end{equation}
Here $a_0 \ll 1$ and $b_0 \gg 1$ are positive constants which will be
chosen later. In \eqref{p0def}, we use a cut-off function $\psi :
(-\infty,1] \to \real$ which satisfies $\psi(y) = 0$ for $y < -1$,
$\psi(-1) = \psi'(-1) = 0$, $0 < \psi'(y) < \frac14(1-\psi(y))$ for
$|y| < 1$, $\psi(1) = 1 - e^{-1/4}$, and $\psi'(1) = \frac14
e^{-1/4}$. For definiteness, we can take
\[
  \psi(y) \,=\, \zeta_2(y+1)^2 - \zeta_3 (y+1)^3 \quad \hbox{for}
  \quad |y| \le 1~,
\]
where $\zeta_2 = \frac34 - \frac78 e^{-1/4} \approx 0.068$ and 
$\zeta_3 = \frac14 - \frac5{16} e^{-1/4} \approx 0.0066$. It 
follows from these definitions that the function $\xi \mapsto 
p_0(\xi,t)$ is piecewise smooth and nondecreasing along rays. 
Moreover, we have
\begin{equation}\label{nablap0}
  \nabla p_0(\xi,t) \,=\, \left\{\begin{array}{lcl}
  \frac{\xi}{2}\,e^{|\xi|^2/4}\Bigl(1 - \tilde\psi(|\xi|^2 - a(t)^2)
  \Bigr) & \hbox{if} & 0 \le |\xi|^2 \le a(t)^2{+}1~, \\[1mm]
  0 & \hbox{if} & a(t)^2{+}1 \le |\xi|^2 \le b(t)^2~, \\[1.5mm]
  \frac{\beta}{4} \frac{\xi}{|\xi|}\,e^{\beta|\xi|/4} & \hbox{if} & 
  |\xi|^2 \ge b(t)^2~,
  \end{array}\right.  
\end{equation}
where $\tilde \psi(y) = \psi(y) + 4\psi'(y)$. The graph of the 
function $p_0(\xi,t)$ is depicted in Fig.~1. 
\figurewithtex 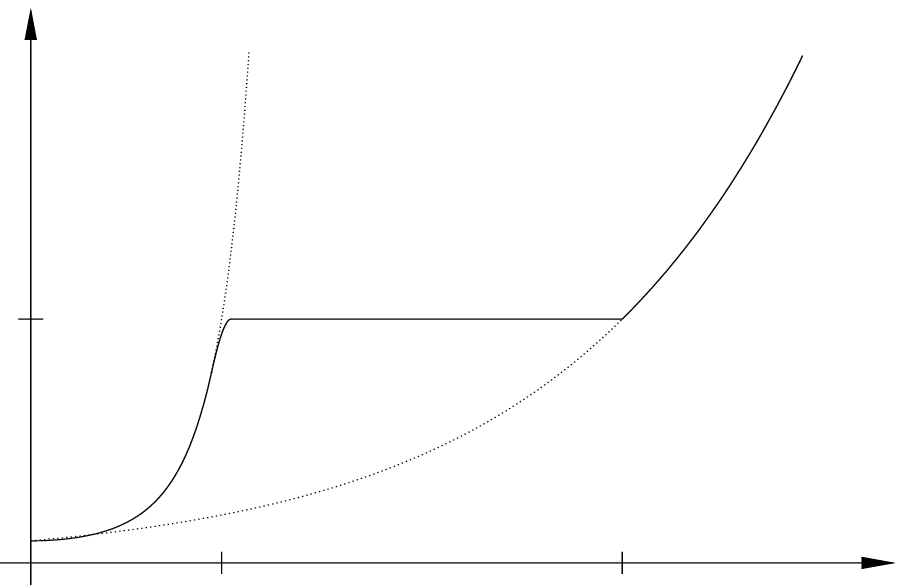 fig.tex 6.000 9.000
{\bf Fig.~1:} The radially symmetric, time-dependent
weight $p_0(\xi,t)$ is represented as a function of $|\xi|$. 
The gradient of $p_0(\xi,t)$ has a jump discontinuity at 
$|\xi| = b(t)$, but is Lipschitz continuous near $|\xi| = a(t)$.\cr

Next, for each $i \in \{1,\dots,N\}$, we define
\begin{equation}\label{pidef}
  p_i(\xi,t) \,=\, p_0(\xi,t) + \Bigl(\frac{\nu t}{d^2}\Bigr)
  q_i(\xi,t)~,
\end{equation}
where the correction $q_i(\xi,t)$ vanishes identically when
$|\xi|^2 \ge a(t)^2 + 1$ and satisfies 
\begin{equation}\label{qidef1}
  v^G(\xi) \cdot \nabla q_i(\xi,t) \,=\, \frac{1}{\pi}
  \sum_{j\neq i}\frac{\alpha_j}{\alpha_i}\frac{d^2}{|z_{ij}(t)|^4}
  (\xi\cdot z_{ij}(t))(\nabla p_0(\xi,t)\cdot z_{ij}(t)^\perp)~,
\end{equation}
for $|\xi|^2 \le a(t)^2 + 1$. Proceeding as in Section~\ref{sec3},
it is easy to obtain the explicit expression
\begin{equation}\label{qidef2}
  q_i(\xi,t) \,=\, -\frac{1}{2}\,\frac{\,\xi\cdot\nabla p_0(\xi,t)}
  {1 - e^{-|\xi|^2/4}}\sum_{j\neq i}\frac{\alpha_j}{\alpha_i} 
  \frac{d^2}{|z_{ij}(t)|^4}\Bigl((\xi\cdot z_{ij}(t))^2 - 
  (\xi\cdot z_{ij}(t)^\perp)^2\Bigr)~,
\end{equation}
which is in fact valid for $|\xi| \le b(t)$, because $\nabla 
p_0(\xi,t)$ vanishes when $a(t)^2 + 1 \le |\xi|^2 \le b(t)^2$. 
It follows from \eqref{qidef2} that
\begin{equation}\label{qibdd}
 (1+|\xi|^2)|q_i(\xi,t)| + |\xi| |\nabla q_i(\xi,t)| + 
 |t\partial_t q_i(\xi,t)| \,\le\, C|\xi|^2(1+|\xi|^2)^2
 p_0(\xi,t)~, 
\end{equation}
whenever $|\xi|^2 \le a(t)^2 + 1$. Using \eqref{pidef} and the
definition \eqref{abdef} of $a(t)$, we deduce from \eqref{qibdd} that
$|p_i(\xi,t) - p_0(\xi,t)| \le C a_0^4 \,p_0(\xi,t)$ for some $C >
0$. Thus, if we choose the constant $a_0 > 0$ sufficiently small, we
see that the weight $p_i(\xi,t)$ is a very small perturbation of
$p_0(\xi,t)$ for all $\xi \in \real^2$ and all $t \in (0,T]$. In
particular, for $i = 0,\dots,N$, we have the uniform bounds
\begin{equation}\label{pbdd}
  \frac12 \,e^{\beta|\xi|/4} \,\le\, p_i(\xi,t) \,\le\, 
  2\,e^{|\xi|^2/4}~, \qquad \xi \in \real^2~, \quad t \in (0,T]~.
\end{equation}

\subsection{Energy estimates}\label{4.2}

Now that we have defined appropriate weights $p_1(\xi,t), \dots, 
p_N(\xi,t)$, it remains to control the evolution of the energy 
functional $E(t)$ introduced in \eqref{Edef}. Here and in what 
follows, we always assume that the parameters in \eqref{abdef} 
satisfy $a_0 \ll 1$, $b_0 \gg 1$, and that the viscosity $\nu > 0$ 
is small enough so that $\nu T \ll d^2$. 

\begin{proposition}\label{energy}
There exist positive constants $\epsilon_0$, $\epsilon_1$ and 
$\kappa_1$, $\kappa_2$, $\kappa_3$, depending only on the ratio $T/T_0$, 
such that, if $a_0 = b_0^{-1} = \epsilon_0$ and $\nu T/d^2 \le 
\epsilon_1$, and if $\tilde w = (\tilde w_1,\dots,\tilde w_N) \in 
C^0((0,T],Y^N)$ is any solution of system \eqref{T0}--\eqref{T4}, 
then the energy $E(t)$ defined in \eqref{Edef} satisfies
\begin{equation}\label{Edecay}
  t E'(t) \,\le\, -\kappa_1 E(t) + \kappa_2\frac{t}{T_0}\Bigl(
  E(t) + E(t)^3\Bigr) + \kappa_3 \Bigl(\frac{\nu t}{d^2}\Bigr)~, 
  \qquad 0 < t \le T~.
\end{equation}
\end{proposition}

\medskip\proof 
Differentiating \eqref{Edef} with respect to time and using
\eqref{T0}--\eqref{T4}, we find
\begin{equation}\label{Eprime}
  tE'(t) \,=\, \sum_{i=1}^N \inttwo \Bigl\{\frac12 \,t \partial_t 
  p_i(\xi,t) |\tilde w_i(\xi,t)|^2 + p_i(\xi,t) \tilde w_i(\xi,t)
  \,t\partial_t \tilde w_i(\xi,t)\Bigr\}\d\xi
  \,=\, \sum_{k=1}^6 \cE_k(t)~,
\end{equation}
where
\begin{align}\nonumber
  \cE_1(t) \,&=\, \sum_{i=1}^N \inttwo \Bigl\{\frac12 \,t \partial_t 
  p_i(\xi,t) |\tilde w_i(\xi,t)|^2 + p_i(\xi,t) \tilde w_i(\xi,t)
  (\cL \tilde w_i(\xi,t) - \tilde w_i(\xi,t))\Bigr\}\d\xi~, \\ \nonumber 
  \cE_2(t) \,&=\, -\sum_{i=1}^N \inttwo p_i(\xi,t)\tilde w_i(\xi,t) 
  \frac{\alpha_i}{\nu} \Bigl(v_i^\app(\xi,t)\cdot \nabla \tilde w_i(\xi,t) 
  + \tilde v_i(\xi,t)\cdot \nabla w_i^\app(\xi,t)\Bigr)\d\xi~, 
  \\ \nonumber
  \cE_3(t) \,&=\, -\sum_{i=1}^N \inttwo p_i(\xi,t)\tilde w_i(\xi,t)
  \sum_{j\neq i} \frac{\alpha_j}{\nu}\left\{v_j^\app \Bigl(\xi + 
  \frac{z_{ij}(t)}{\sqrt{\nu t}}\,,\,t\Bigr) - v^G\Bigl(\frac{z_{ij}(t)}
  {\sqrt{\nu t}}\Bigr)\right\}\cdot \nabla \tilde w_i(\xi,t)\d\xi~,
  \\ \nonumber
  \cE_4(t) \,&=\, -\sum_{i=1}^N \inttwo p_i(\xi,t)\tilde w_i(\xi,t)
  \sum_{j\neq i} \frac{\alpha_j}{\nu}\,\tilde v_j \Bigl(\xi + 
  \frac{z_{ij}(t)}{\sqrt{\nu t}}\,,\,t\Bigr)\cdot \nabla w_i^\app(\xi,t)
  \d\xi~, \\ \nonumber
  \cE_5(t) \,&=\, -\sum_{i=1}^N \inttwo p_i(\xi,t)\tilde w_i(\xi,t)
  \sum_{j=1}^N \frac{\alpha_j t}{d^2}\,\tilde v_j\Bigl(\xi + 
  \frac{z_{ij}(t)}{\sqrt{\nu t}}\,,\,t\Bigr)\cdot \nabla \tilde w_i(\xi,t)
  \d\xi~, \\ \nonumber  
  \cE_6(t) \,&=\, -\sum_{i=1}^N \inttwo p_i(\xi,t)\tilde w_i(\xi,t)
  \tilde R_i(\xi,t)\d\xi~.
\end{align}
The general strategy is to control the contributions of $\cE_2(t), 
\dots, \cE_6(t)$ using the negative terms contained in $\cE_1(t)$. 
We shall proceed in six steps, each one being devoted to the detailed 
analysis of a specific term. Very often, we will have to consider 
separately the three cases $|\xi|^2 \le a(t)^2 + 1$, $a(t)^2 + 1 \le 
|\xi|^2 \le b(t)^2$, and $|\xi| \ge b(t)$, see \eqref{p0def}. 
For simplicity, these three domains in $\real^2$ will be denoted by 
``region I'', ``region II'', and ``region III'', respectively. 

\medskip\noindent{\bf Step 1:}  {\sl Diffusive terms} \\
Our goal is to show that there exists a constant $\kappa > 0$ such that
\begin{equation}\label{EE1}
  \cE_1(t) \,\le\, -\kappa \sum_{i=1}^N \Bigl(\cI_i(t) + 
  \cJ_i(t) + \cK_i(t)\Bigr)~,
\end{equation}
where 
\begin{align}\nonumber
  \cI_i(t) \,&=\, \inttwo p_0(\xi,t) |\nabla \tilde w_i(\xi,t)|^2
  \d\xi~, \qquad
  \cJ_i(t) \,=\, \inttwo p_0(\xi,t)\,\chi(\xi,t)|\tilde w_i(\xi,t)|^2
  \d\xi~, \\ \label{cKidef}
  \cK_i(t) \,&=\, \inttwo p_0(\xi,t)|\tilde w_i(\xi,t)|^2\d\xi~.
\end{align}
Here $\chi(\xi,t)$ is the continuous, radially symmetric function
defined by
\[
  \chi(\xi,t) \,=\, \left\{\begin{array}{lcl}
  |\xi|^2 & \hbox{if} & 0 \le |\xi| \le a(t)~, \\
  a(t)^2 & \hbox{if} & a(t) \le |\xi| \le b(t)~, \\
  \beta|\xi| & \hbox{if} & |\xi| \ge b(t)~.
  \end{array}\right.  
\]
To prove \eqref{EE1}, we start from the identity
\[
  \inttwo p_i \tilde w_i (\cL\tilde w_i - \tilde w_i)\d\xi 
  \,=\, -\inttwo \Bigl(p_i |\nabla\tilde w_i|^2 + \tilde w_i
  (\nabla p_i \cdot \nabla\tilde w_i) + \frac14 (\xi\cdot\nabla 
  p_i)|\tilde w_i|^2 + \frac12 p_i |\tilde w_i|^2\Bigr)\d\xi~,   
\]
which is easily obtained using \eqref{cLdef} and integrating by parts. 
Since
\[
  \Bigl|\inttwo \tilde w_i(\nabla p_i \cdot \nabla\tilde w_i)\d\xi
  \Bigr| \,\le\, \frac34 \inttwo p_i |\nabla\tilde w_i|^2\d\xi
  + \frac13 \inttwo \frac{|\nabla p_i|^2}{p_i}\,|\tilde w_i|^2
  \d\xi~,
\]
we see that $\cE_1(t) \le \tilde \cE_1(t)$, where
\[
  \tilde \cE_1(t) \,=\, - \sum_{i=1}^N \inttwo \Bigl(\frac14 p_i(\xi,t) 
  |\nabla\tilde w_i(\xi,t)|^2 + p_i(\xi,t) \tilde \chi(\xi,t) 
  |\tilde w_i(\xi,t)|^2 + \frac12 p_i(\xi,t) |\tilde w_i(\xi,t)|^2
  \Bigr)\d\xi~,
\]
and
\[
  \tilde \chi(\xi,t) \,=\, \frac{\xi\cdot\nabla p_i(\xi,t)}{4 p_i(\xi,t)}
  -\frac{|\nabla p_i(\xi,t)|^2}{3 p_i(\xi,t)^2} - 
  \frac{t\partial_t p_i(\xi,t)}{2 p_i(\xi,t)}~.
\]
We already observed that $p_i(\xi,t) \ge \frac12 p_0(\xi,t)$ if $a_0$ 
is sufficiently small. So, to prove \eqref{EE1}, it remains to show 
that $\tilde \chi(\xi,t) \ge C \chi(\xi,t)$ for some $C > 0$. This 
is easily verified in regions II and III, because $p_i(\xi,t) = 
p_0(\xi,t)$ in these regions. From \eqref{p0def}, we find by 
a direct calculation
\[
  \tilde \chi(\xi,t) \,=\, \left\{\begin{array}{lcl}
  {\DS \frac{a(t)^2}{16}} & \hbox{ if } & a(t)^2 + 1 \le |\xi|^2 \le 
  b(t)^2~, \\[2mm]
  {\DS \frac{\beta}{16}|\xi| - \frac{\beta^2}{48}} & \hbox{ if } & 
  |\xi| \ge b(t)~. \end{array}\right.
\]
In region I, we first compute the contributions of the radially 
symmetric weight $p_0(\xi,t)$ to the function $\tilde \chi(\xi,t)$. 
Using \eqref{nablap0}, we obtain
\[
  \tilde \chi_0(\xi,t) \,=\, \frac{|\xi|^2}{8} \frac{1-\tilde 
  \psi(|\xi|^2{-}a^2)}{1-\psi(|\xi|^2{-}a^2)} -\frac{|\xi|^2}{12}
  \Bigl(\frac{1 -\tilde \psi(|\xi|^2{-}a^2)}{1 -\psi(|\xi|^2{-}a^2)}
  \Bigr)^2 + \frac{a(t)^2}{4}\frac{\psi'(|\xi|^2{-}a^2)}{1 - 
  \psi(|\xi|^2{-}a^2)}~.
\]
In particular, we see that $\tilde \chi_0(\xi,t) = |\xi|^2/(24)$
if $|\xi|^2 \le a(t)^2-1$. If $y = |\xi|^2 - a(t)^2 \in [-1,1]$, 
we use the elementary bounds
\[
  0 \le \frac{1 - \tilde \psi(y)}{1 - \psi(y)} \,\le\, 1~, 
  \quad \hbox{and}\quad 
  \frac{\frac{1}{24}(1-\tilde \psi(y)) + \frac14\psi'(y)}{1-\psi(y)} 
  \,\ge\, \delta \,>\, 0~,
\]
which follow from the definition of $\psi$, and we deduce that $\tilde
\chi_0(\xi,t) \ge \delta\min\{|\xi|^2,a(t)^2\}$ for some $\delta > 0$.
Summarizing, we have shown that there exists $C_0 > 0$ such that 
$\tilde \chi_0(\xi,t) \ge C_0 |\xi|^2$ in region I. On the other hand, 
using \eqref{qibdd}, it is straightforward to verify that $\tilde 
\chi(\xi,t) \ge \tilde \chi_0(\xi,t) - C_1 a_0^4|\xi|^2$ when 
$|\xi|^2 \le a(t)^2+1$. Thus, if the constant $a_0$ in \eqref{abdef} is 
sufficiently small, there exists $C_2 > 0$ such that $\tilde \chi(\xi,t) 
\ge C_2\chi(\xi,t)$ for all $\xi \in \real^2$ and all $t \in (0,T]$. 
This concludes the proof of \eqref{EE1}.

\medskip\noindent{\bf Step 2:} {\sl Self-interaction terms}\\
We next consider the term $\cE_2(t)$ in \eqref{Eprime}. To simplify 
the notations, we rewrite \eqref{wvapp} in the form
\begin{equation}\label{wvaux}
  w_i^\app(\xi,t) \,=\, G(\xi) + \Bigl(\frac{\nu t}{d^2}\Bigr)
  \cF_i(\xi,t)~, \qquad 
  v_i^\app(\xi,t) \,=\,  v^G(\xi) + \Bigl(\frac{\nu t}{d^2}\Bigr)
  v^{\cF_i}(\xi,t)~,
\end{equation} 
where $\cF_i(\xi,t) = F_i(\xi,t) + (\nu t/d^2)^{1/2}H_i(\xi,t) + 
(\nu t/d^2)K_i(\xi,t)$, and $v^{\cF_i}$ is the velocity field
obtained from $\cF_i$ via the Biot-Savart law \eqref{BS}. We thus
have $\cE_2(t) = \Omega_1(t) + \dots + \Omega_4(t)$, where
\begin{align*}
  \Omega_1(t) \,&=\, -\sum_{i=1}^N \frac{\alpha_i}{\nu} \inttwo 
  p_i(\xi,t)\tilde w_i(\xi,t)(v^G(\xi) \cdot \nabla \tilde w_i(\xi,t))
  \d\xi~,\\ 
  \Omega_2(t) \,&=\, -\sum_{i=1}^N \frac{\alpha_i}{\nu} \inttwo 
  p_0(\xi,t)\tilde w_i(\xi,t)(\tilde v_i(\xi,t) \cdot \nabla G(\xi))
  \d\xi~,\\
  \Omega_3(t) \,&=\, -\sum_{i=1}^N \frac{\alpha_i t}{d^2}
  \inttwo q_i(\xi,t)\tilde w_i(\xi,t)(\tilde v_i(\xi,t) \cdot \nabla 
  G(\xi))\d\xi~,\\
  \Omega_4(t)\,&=\, -\sum_{i=1}^N \frac{\alpha_i t}{d^2} \inttwo 
  p_i(\xi,t)\tilde w_i(\xi,t)\Bigl(v^{\cF_i}(\xi,t)\cdot \nabla 
  \tilde w_i(\xi,t) + \tilde v_i(\xi,t)\cdot \nabla \cF_i(\xi,t)
  \Bigr)\d\xi~.
\end{align*}
To prove that $\Omega_1(t)$ has a finite limit as $\nu \to 0$, 
we integrate by parts and use the fact that $v^G\cdot\nabla p_0 = 0$
because the weight $p_0$ is radially symmetric. In view of 
\eqref{pidef}, we find
\[
  \Omega_1(t) \,=\, \sum_{i=1}^N \frac{\alpha_i}{2\nu} \inttwo 
  |\tilde w_i(\xi,t)|^2 (v^G(\xi)\cdot \nabla p_i(\xi,t))\d\xi 
  \,=\, \sum_{i=1}^N \frac{\alpha_i t}{2d^2} \inttwo |\tilde 
  w_i(\xi,t)|^2 (v^G(\xi) \cdot \nabla q_i(\xi,t))\d\xi~. 
\]
This term cannot be controlled by \eqref{EE1}, unless $t \ll T_0$, 
but it will exactly compensate for another term coming from $\cE_3(t)$. 
As was already explained, this is precisely the reason why the 
correction $(\nu t/d^2)q_i(\xi,t)$ was added to the weight 
$p_0(\xi,t)$. To treat $\Omega_2$, we use the fact that the 
linear operator $\tilde w_i \mapsto \tilde v_i\cdot\nabla G$ 
(where $\tilde v_i$ is obtained from $\tilde w_i$ via the Biot-Savart 
law) is skew-symmetric in the function space $Y$ defined by \eqref{Ydef}, 
see \cite{GW05}. We thus have
\begin{align*}
  |\Omega_2(t)| \,&=\, \biggl|\sum_{i=1}^N \frac{\alpha_i}{\nu} 
  \inttwo \Bigl(e^{|\xi|^2/4} - p_0(\xi,t)\Bigr)\tilde w_i(\xi,t)
  (\tilde v_i(\xi,t) \cdot \nabla G(\xi))\d\xi\biggr| \\
  \,&\le\, \sum_{i=1}^N \frac{|\alpha_i|}{\nu} 
  \int_{|\xi|^2 \ge a(t)^2-1} e^{|\xi|^2/4}|\tilde w_i(\xi,t)| 
  \,|\tilde v_i(\xi,t)|\,|\nabla G(\xi)| \d\xi \\
  \,&\le\, C\sum_{i=1}^N \frac{|\alpha_i|}{\nu}
  \int_{|\xi|^2 \ge a(t)^2-1} e^{\beta|\xi|/8}|\tilde w_i(\xi,t)| 
  \,|\tilde v_i(\xi,t)|\,|\xi| e^{-\beta|\xi|/8}\d\xi~.
\end{align*}
In the last line, we have used the definition \eqref{Gdef} of $G$.  To
estimate the last integral, we apply the trilinear H\"older inequality
with exponents $2$, $4$, $4$. We have $\|e^{\beta|\xi|/8} \tilde
w_i\|_{L^2} \le C\cK_i^{1/2}$ by \eqref{pbdd}, \eqref{cKidef}, and
$\|\tilde v_i\|_{L^4} \le C\|\tilde w_i\|_{L^{4/3}} \le C\cK_i^{1/2}$ 
(see e.g. \cite{GW02}). Moreover, by direct calculation, 
\[
  \Bigl(\int_{|\xi|^2 \ge a(t)^2-1}|\xi|^4\,e^{-\beta|\xi|/2}\d\xi
  \Bigr)^{1/4} \,\le\, C\,a(t)^{3/2}\,\frac{e^{-\beta a(t)/8}}
  {(\beta a(t))^{1/4}}~, 
\]
provided that $\beta a(t) \ge 1$. Thus, if the constant $\epsilon_1$ 
in Proposition~\ref{energy} is sufficiently small (depending on 
$\epsilon_0$ and $T/T_0$), we find
\[
  |\Omega_2(t)| \,\le\, C \sum_{i=1}^N \frac{|\alpha_i|}{\nu} 
  \,\cK_i(t)\,a(t)^{3/2}\,e^{-\beta a(t)/8} \,\le\, \epsilon 
  \sum_{i=1}^N \cK_i(t)~. 
\]
Here and in what follows, $\epsilon$ denotes a positive constant
which can be made arbitrarily small by an appropriate choice 
of $\epsilon_0$ and $\epsilon_1$. Using similar estimates, it is 
also easy to bound the regular terms $\Omega_3$ and $\Omega_4$. 
We find
\begin{align*}
  |\Omega_3(t)| \,&\le\, C\sum_{i=1}^N \frac{|\alpha_i| t}{d^2}
  \inttwo p_0(\xi,t)|\xi|^2(1+|\xi|^2)|\tilde w_i(\xi,t)|
  \,|\tilde v_i(\xi,t)|\,|\nabla G(\xi)| \d\xi \,\le\, 
  C\frac{t}{T_0}\sum_{i=1}^N \cK_i(t)~, \\
  |\Omega_4(t)| \,&\le\, C\sum_{i=1}^N \frac{|\alpha_i| t}{d^2} 
  \inttwo p_0(\xi,t)|\tilde w_i(\xi,t)|\Bigl(|v^{\cF_i}(\xi,t)| 
  \,|\nabla \tilde w_i(\xi,t)| + |\tilde v_i(\xi,t)|\,|\nabla\cF_i(\xi,t)|
  \Bigr)\d\xi \\
  \,&\le\, C\sum_{i=1}^N \frac{|\alpha_i| t}{d^2}\Bigl(\cK_i(t)^{1/2}
  \cI_i(t)^{1/2} + \cK_i(t)\Bigr) \,\le\, \epsilon \sum_{i=1}^N
  \cI_i(t) + C\frac{t}{T_0}\sum_{i=1}^N\cK_i(t)~.
\end{align*}

\medskip\noindent{\bf Step 3:}  {\sl Advection terms}\\
Using \eqref{wvaux} and integrating by parts,  we write $\cE_3(t) 
= \Psi_1(t) + \Psi_2(t) + \Psi_3(t)$, where
\begin{align*}
  \Psi_1(t) \,&=\, \sum_{i=1}^N \inttwo |\tilde w_i(\xi,t)|^2 
  \,\nabla p_0(\xi,t) \cdot \sum_{j\neq i} \frac{\alpha_j}{2\nu}
  \left\{v^G\Bigl(\xi + \frac{z_{ij}(t)}{\sqrt{\nu t}}\Bigr) 
  - v^G\Bigl(\frac{z_{ij}(t)}{\sqrt{\nu t}}\Bigr)\right\}\d\xi~, \\
  \Psi_2(t) \,&=\, \sum_{i=1}^N \inttwo |\tilde w_i(\xi,t)|^2 
  \,\nabla q_i(\xi,t) \cdot \sum_{j\neq i} \frac{\alpha_j t}{2d^2}
  \left\{v^G\Bigl(\xi + \frac{z_{ij}(t)}{\sqrt{\nu t}}\Bigr) 
  - v^G\Bigl(\frac{z_{ij}(t)}{\sqrt{\nu t}}\Bigr)\right\}\d\xi~, \\
  \Psi_3(t) \,&=\, \sum_{i=1}^N \inttwo |\tilde w_i(\xi,t)|^2 
  \,\nabla p_i(\xi,t) \cdot \sum_{j\neq i} \frac{\alpha_j t}{2d^2} 
  \,v^{\cF_j}\Bigl(\xi + \frac{z_{ij}(t)}{\sqrt{\nu t}}\,,\,t\Bigr) 
  \d\xi~.
\end{align*}
By construction, all three integrands vanish for $a(t)^2 + 1 \le 
|\xi|^2 \le b(t)^2$, so we need only consider regions I and III. 
The most important term is the contribution of region I to 
$\Psi_1(t)$, which we denote by $\Psi_1^I(t)$. Using 
\eqref{Deltaidef} together with the identity
\[
  v^G(\xi)\cdot\nabla q_i(\xi,t) \,=\, -\nabla p_0(\xi,t)\cdot
  D_i(\xi,t)~, \qquad |\xi|^2 \le a(t)^2+1~,
\]
which follows immediately from the definitions \eqref{Didef}, 
\eqref{qidef1}, we find
\[
  \Psi_1^I(t) \,=\, -\Omega_1(t) + \frac12 \sum_{i=1}^N \int_{|\xi|^2 \le 
  a(t)^2+1} |\tilde w_i(\xi,t)|^2 \,\nabla p_0(\xi,t) \cdot 
  \Delta_i(\xi,t) \d \xi~.
\]
Thus, using the estimates on $\Delta_i(\xi,t)$ which follow from
\eqref{Deltaiexp}, \eqref{Wexp}, we obtain
\begin{align*}
  |\Psi_1^I(t) + \Omega_1(t)| \,&\le\, C\sum_{i=1}^N \sum_{j\neq i} 
  \frac{|\alpha_j|}{\nu}\int_{|\xi|^2 \le a(t)^2+1} p_0(\xi,t) 
  |\tilde w_i(\xi,t)|^2 |\xi|^3 \Bigl(\frac{\nu t}{d^2}\Bigr)^{3/2}\d\xi \\
  \,&\le\, C \frac{|\alpha|t}{d^2}\,a(t)\Bigl(\frac{\nu t}{d^2}\Bigr)^{1/2}
  \sum_{i=1}^N \cJ_i(t) \,\le\, \epsilon \sum_{i=1}^N \cJ_i(t)~.
\end{align*}
On the other hand, the contribution of region III to $\Psi_1(t)$ 
can be estimated as follows
\begin{align*}
  |\Psi_1^{III}(t)| \,&\le\, C\sum_{i=1}^N \sum_{j\neq i} 
  \frac{|\alpha_j|}{\nu}\int_{|\xi| \ge b(t)} \beta p_0(\xi,t) 
  |\tilde w_i(\xi,t)|^2 \Bigl(\frac{\nu t}{d^2}\Bigr)^{1/2}\d\xi \\
  \,&\le\, C \frac{|\alpha|t}{d^2}\,\frac{1}{b_0}\sum_{i=1}^N
  \int_{|\xi| \ge b(t)} \beta |\xi| p_0(\xi,t) |\tilde w_i(\xi,t)|^2 
  \d\xi \,\le\, \epsilon \sum_{i=1}^N \cJ_i(t)~,
\end{align*}
if the parameter $b_0 > 0$ is chosen sufficiently large. Similarly, 
using \eqref{qibdd}, we can bound $\Psi_2(t)$ in the following 
way
\begin{align*}
  |\Psi_2(t)| \,&\le\, C \sum_{i=1}^N \sum_{j\neq i}\frac{|\alpha_j| t}{d^2}
  \int_{|\xi|^2 \le a(t)^2+1} p_0(\xi,t) |\xi|^2(1+|\xi|^2)^2
  |\tilde w_i(\xi,t)|^2 \Bigl(\frac{\nu t}{d^2}\Bigr)\d\xi \\
  \,&\le\, C\frac{|\alpha|t}{d^2}\,a_0^4 \sum_{i=1}^N 
  \int_{|\xi|^2 \le a(t)^2+1} p_0(\xi,t) |\xi|^2|\tilde w_i(\xi,t)|^2 
  \d\xi \,\le\, \epsilon \sum_{i=1}^N \cJ_i(t)~,
\end{align*}
if $a_0 > 0$ is sufficiently small. Finally, to bound $\Psi_3(t)$, we
recall that $|v^{\cF_i}(\xi,t)| \le C(1+|\xi|^2)^{-3/2}$, see
\eqref{Fibdd}, \eqref{Hibdd}, \eqref{Kibdd}. Since $|\xi + 
\frac{z_{ij}}{\sqrt{\nu t}}| \ge \frac{d}{2\sqrt{\nu t}}$ in regions
I and III when $i \neq j$, we find
\[
  |\Psi_3(t)| \,\le\, C \sum_{i=1}^N \sum_{j\neq i} \frac{|\alpha_j| t}
  {d^2}\Bigl(\frac{\nu t}{d^2}\Bigr)^{3/2} \inttwo |\nabla p_0(\xi,t)|
  |\tilde w_i(\xi,t)|^2 \d\xi \,\le\, \epsilon \sum_{i=1}^N \cK_i(t)~.
\]

\medskip\noindent{\bf Step 4:}  {\sl Cross-interaction terms}\\
To bound $\cE_4(t)$, we need a good estimate on the velocity 
field $\tilde v_i(\xi,t)$. We claim that there exists a constant
$C > 0$ such that
\begin{equation}\label{vjest}
  \|(1+|\xi|^2)\tilde v_i(\xi,t)\|_{L^\infty}^2 \,\le\, 
  C\Bigl(\cI_i(t) + \cK_i(t)\Bigr)~, \qquad i \in \{1,\dots,N\}~. 
\end{equation}
Indeed, since $\inttwo \tilde w_i(\xi,t)\d\xi \,=\, 0$ for all $t \in 
(0,T]$, it follows from \cite[Proposition~B.1]{GW02} that
\[
  \|(1+|\xi|^2)\tilde v_i\|_{L^\infty}^2 \,\le\, C\Bigl(\|(1+|\xi|^2) 
  \tilde w_i\|_{L^4} + \|(1+|\xi|^2)\tilde w_i\|_{L^{4/3}}\Bigr)~.
\]
As $H^1(\real^2) \hookrightarrow L^4(\real^2)$ and $(1+|\xi|^2) \le 
Cp_0(\xi,t)$, we have $\|(1+|\xi|^2) \tilde w_i\|_{L^4} \le 
C(\cI_i^{1/2} + \cK_i^{1/2})$. In addition, using H\"older's inequality,
we find $\|(1+|\xi|^2)\tilde w_i\|_{L^{4/3}} \le C\cK_i^{1/2}$. 
This proves \eqref{vjest}. 

The main contribution to $\cE_4(t)$ comes from region I. Since 
$|\xi + \frac{z_{ij}}{\sqrt{\nu t}}| \ge \frac{d}{2\sqrt{\nu t}}$ 
when $i \neq j$ and $|\xi|^2 \le a(t)^2+1$, we obtain, using
\eqref{vjest},
\begin{align*}
  |\cE_4^I(t)| \,&\le\, C\sum_{i=1}^N \sum_{j\neq i} 
  \frac{|\alpha_j| t}{d^2}\Bigl(\cI_j(t) + \cK_j(t)\Bigr)^{1/2}
  \int_{|\xi|^2 \le a(t)^2+1} p_0(\xi,t)|\tilde w_i(\xi,t)|
  |\nabla w_i^\app(\xi,t)|\d\xi \\
  \,&\le\, C\frac{|\alpha| t}{d^2}\sum_{j=1}^N\Bigl(\cI_j(t) + 
  \cK_j(t)\Bigr)^{1/2}\sum_{i=1}^N\cK_i(t)^{1/2} \,\le\, \epsilon
  \sum_{i=1}^N\cI_i(t) + C\frac{t}{T_0}\sum_{i=1}^N\cK_i(t)~.
\end{align*}
In regions II and III, the quantity $|\nabla w_i^\app(\xi,t)|$ is 
bounded by $Ce^{-\gamma |\xi|^2/4}$ for any $\gamma < 1$. Choosing 
$\gamma > 1/2$ and proceeding as in the second step, we easily find
\[
  |\cE_4^{II}(t)| + |\cE_4^{III}(t)| \,\le\, C\sum_{i=1}^N \sum_{j\neq i} 
  \frac{|\alpha_j|}{\nu}\,\cK_i(t)^{1/2}\cK_j(t)^{1/2} e^{-(\gamma-\frac12)
  a(t)^2/4} \,\le\, \epsilon \sum_{i=1}^N\cK_i(t)~.
\]

\medskip\noindent{\bf Step 5:}  {\sl Nonlinear terms}\\
Instead of \eqref{vjest}, we use here the simpler inequality
\[
  \|\tilde v_i\|_{L^\infty}^2 \,\le\, C \|\tilde w_i\|_{L^4} 
  \|\tilde w_i\|_{L^{4/3}} \,\le\, C\Bigl(\cI_i + \cK_i\Bigr)^{1/2}
  \cK_i^{1/2}~, 
\]
which follows from \cite[Lemma 2.1]{GW02}. Applying H\"older's 
inequality, we thus obtain
\begin{align*}
  |\cE_5(t)| \,&\le\, C\sum_{i,j=1}^N \frac{|\alpha_j| t}{d^2}
  \cK_i(t)^{1/2}\Bigl(\cI_j(t) + \cK_j(t)\Bigr)^{1/4}\cK_j(t)^{1/4}
  \cI_i(t)^{1/2} \\
  \,&\le\, \epsilon \sum_{i=1}^N \cI_i(t) + C\frac{t}{T_0}
  \sum_{i=1}^N \Bigl(\cK_i(t)^2 + \cK_i(t)^3\Bigr)~.
\end{align*}

\medskip\noindent{\bf Step 6:}  {\sl Remainder terms}\\
Finally, using H\"older's inequality and estimate \eqref{tildeRi}, 
we find
\[
  |\cE_6(t)| \,\le\, C\sum_{i=1}^N \Bigl(\frac{\nu t}{d^2}\Bigr)^{1/2} 
  \cK_i(t)^{1/2} \,\le\, \epsilon \sum_{i=1}^N \cK_i(t) + 
  C\Bigl(\frac{\nu t}{d^2}\Bigr)~.
\]
Collecting the estimates established in Steps 1--6 and using the fact 
that $E(t) \approx \frac12\sum_{i=1}^N \cK_i(t)$, we easily obtain 
\eqref{Edecay} provided that $\epsilon > 0$ is sufficiently small. 
As was already explained, this last condition is easy to fulfill 
if we choose the constants $\epsilon_0$, $\epsilon_1$ appropriately.  
This concludes the proof of Proposition~\ref{energy}. \QED

\subsection{End of the proof of Theorem~\ref{thm4}}\label{4.3}

It is now quite easy to conclude the proof of Theorem~\ref{thm4}.
If the solution $\omega^\nu(x,t)$ of \eqref{V} with initial data 
\eqref{mudef} is decomposed as in \eqref{omdecomp}, we know from
Lemma~\ref{lemdec} that the rescaled vorticity profiles 
$w_i(\xi,t) \equiv w_i^\nu(\xi,t)$ defined by \eqref{wvdef}
satisfy $w_i(\cdot,t) \in Y$ for all $t \in (0,T]$, see 
\eqref{omegaibdd}. Standard parabolic estimates then imply 
that $w = (w_1,\dots,w_N) \in C^0((0,T],Y^N)$ is a solution 
of system \eqref{wieq}, where the vortex positions $z_i(t) \equiv 
z_i^\nu(t)$ are given by \eqref{PW3}. Moreover, we know from 
\cite[Proposition~4.5]{GG05} that $w_i(\cdot,t) \to G$ in 
$Y$ as $t \to 0$, for all $i \in \{1,\dots,N\}$. In fact, in 
Section~6.3 of \cite{GG05}, this convergence is established 
in a polynomially weighted space only, but the proof also works
(and is in fact simpler) in the Gaussian space $Y$. Using the 
approximate solution $w_i^\app(\xi,t)$ of \eqref{wieq} constructed
in Section~\ref{sec3}, we can even obtain the following improved 
estimate for short times:

\begin{lemma}\label{smallt}
For any fixed $\nu > 0$, one has $\|w_i(\cdot,t) - w_i^\app(\cdot,t)
\|_Y = \cO(t^{3/2})$ as $t \to 0$, for all $i \in \{1,\dots,N\}$. 
\end{lemma}

\noindent{\bf Proof.} Let $\hat w_i(\xi,t) = w_i(\xi,t) - 
w_i^\app(\xi,t)$. Then $\hat w_i$ satisfies a system which is 
similar to \eqref{T0}--\eqref{T4}, except that the first line
\eqref{T0} is replaced by $t\partial_t \hat w_i(\xi,t) - (\cL \hat 
w_i)(\xi,t)$, and in the last line \eqref{T4} the left-hand side 
becomes
\begin{equation}\label{hatnl}
   \sum_{j=1}^N \frac{\alpha_j}{\nu}\,\hat v_j
   \Bigl(\xi + \frac{z_{ij}(t)}{\sqrt{\nu t}}\,,\,t\Bigr)
   \cdot \nabla \hat w_i(\xi,t) \,+\, R_i^{(3)}(\xi,t)~,
\end{equation}
where $R_i^{(3)}(\xi,t)$ satisfies \eqref{Ri3}. Note that the nonlinear
terms in \eqref{hatnl} are singular as $\nu \to 0$, but this is not
a problem here because $\nu > 0$ is fixed. To estimate $\hat w_i(\cdot,t)$ 
in the space $Y$, we use the energy functional
\[
  \hat E(t) \,=\, \frac12 \inttwo e^{|\xi|^2/4} \Bigl(
  |\hat w_1(\xi,t)|^2 + \dots + |\hat w_N(\xi,t)|^2\Bigr)\d\xi~. 
\]
Repeating the proof of Proposition~\ref{energy}, with substantial 
simplifications, we obtain for small times a differential inequality 
of the form
\begin{equation}\label{hatEdecay}
  t \hat E'(t) \,\le\, -\eta_1 \hat E(t) + \eta_2 \hat E(t)^3  
  + \eta_3 \Bigl(\frac{\nu t}{d^2}\Bigr)^3~, 
\end{equation}
where the positive constants $\eta_1, \eta_2, \eta_3$ may depend 
on $\nu$. In the derivation of \eqref{hatEdecay}, the only new 
ingredient is the estimate
\[
  \inttwo e^{|\xi|^2/4}\,\hat w_i \cL \hat w_i\d\xi \,\le\, 
  -\kappa \inttwo e^{|\xi|^2/4}\Bigl(|\nabla \hat w_i|^2 + 
  |\xi|^2 |\hat w_i|^2 + |\hat w_i|^2\Bigr)\d\xi~, 
\]
which holds for some $\kappa > 0$ because the self-adjoint 
operator $\cL$ is strictly negative in the subspace of functions
with zero mean, see \cite{GW02}. Since we already know that 
$\hat E(t) \to 0$ as $t \to 0$, inequality \eqref{hatEdecay} 
implies that $\hat E'(t) \le Ct^2$ for $t > 0$ sufficiently 
small, hence $\hat E(t) = \cO(t^3)$ as $t \to 0$. This is the 
desired result. \QED

\medskip
We now consider the energy functional $E(t)$ defined by \eqref{Edef}. 
Since $\hat w_i(\xi,t) = (\nu t/d^2)\tilde w_i(\xi,t)$, and since 
the weights $p_i(\xi,t)$ satisfy \eqref{pbdd}, it follows from 
Lemma~\ref{smallt} that $E(t) \to 0$ as $t \to 0$, for any fixed 
$\nu > 0$. As long as $E(t) \le 1$, we have by Proposition~\ref{energy}
\[
  tE'(t) \,\le\, 2\kappa_2 \,\frac{t}{T_0}E(t) + \kappa_3 
  \Bigl(\frac{\nu t}{d^2}\Bigr)~,
\]
hence 
\begin{equation}\label{Eint}
  E(t) \,\le\, \kappa_3 \int_0^t e^{2\kappa_2s/T_0} 
  \Bigl(\frac{\nu}{d^2}\Bigr)\d s \,\le\, 
  \kappa_3\,e^{2\kappa_2 t/T_0}\Bigl(\frac{\nu t}{d^2}\Bigr)~.
\end{equation}
If we choose $\nu > 0$ sufficiently small so that
\[
  \kappa_3\,e^{2\kappa_2 T/T_0}\Bigl(\frac{\nu T}{d^2}\Bigr) \,\le\, 1~,
\]
we see that \eqref{Eint} holds for all $t \in (0,T]$. Since 
$p_i(\xi,t) \ge Ce^{\beta|\xi|/4}$ for all $\xi \in \real^2$ and
all $t \in (0,T]$, we obtain \eqref{Econc}, and \eqref{strongconv2} 
follows. The proof of Theorem~\ref{thm4} is thus complete. 
\QED


\section{Appendix}\label{app}

\noindent{\bf Proof of Lemma~\ref{lemdec}.} Let $U : \real^2 \times 
(0,\infty) \to \real^2$ be a smooth, divergence-free vector field, 
and fix $\nu > 0$. We assume that
\begin{equation}\label{Ubdd}
  \frac{1}{\nu}\sup_{t>0}\Bigl((\nu t)^{1/2}\|U(\cdot,t)\|_{L^\infty} + 
  \|\Omega(\cdot,t)\|_{L^1}\Bigr) \,=\, C_0 \,<\,\infty~,
\end{equation}
where $\Omega = \partial_1 U_2  - \partial_2 U_1$. Then any solution 
of the linear equation $\partial_t \omega + (U\cdot\nabla)\omega = 
\nu \Delta \omega$ can be represented as
\begin{equation}\label{omegrep}
  \omega(x,t) \,=\, \inttwo \Gamma^\nu_U(x,t;y,s)\omega(y,s)\d y~, 
  \quad x \in \real^2~, \quad t > s > 0~,
\end{equation}
where the {\em fundamental solution} $\Gamma^\nu_U(x,t;y,s)$ has 
the following properties:

\noindent{\bf 1.} For any $\beta \in (0,1)$, there exists $C_1 = 
C_1(\beta,C_0) > 0$ such that
\begin{equation}\label{Gammup}
  0 \,<\, \Gamma^\nu_U(x,t;y,s) \,\le\, \frac{C_1}{\nu(t-s)}
  \,\exp\Bigl(-\beta\,\frac{|x-y|^2}{4\nu (t-s)}\Bigr)~, 
\end{equation}
for all $x,y \in \real^2$ and all $t > s > 0$. This very precise
upper bound is due to E.~Carlen and M.~Loss \cite{CL96}. 

\noindent{\bf 2.} There exists $\gamma = \gamma(C_0) > 0$ and, 
for any $\delta > 0$, there exists $C_2 = C_2(\delta,C_0) > 0$ 
such that
\begin{equation}\label{Gammcont}
  |\Gamma^\nu_U(x,t;y,s) - \Gamma^\nu_U(x',t';y',s')| \,\le\, 
  C_2\Bigl(|x-x'|^\gamma + |t-t'|^{\gamma/2} + |y-y'|^\gamma + 
  |s-s'|^{\gamma/2}\Bigr)~,
\end{equation}
whenever $t-s \ge \delta$ and $t'-s' \ge \delta$. This H\"older
continuity property, which is due to H.~Osada \cite{Os87}, implies
in particular that $\Gamma^\nu_U(x,t;y,s)$ can be continuously 
extended up to $s = 0$, and that this extension (which is still 
denoted by $\Gamma^\nu_U$) satisfies \eqref{Gammup} and 
\eqref{Gammcont} with $s = 0$.

\noindent{\bf 3.} For all $x,y \in \real^2$ and all $t > s > 0$, 
we have
\begin{equation}\label{Gammint}
  \inttwo \Gamma^\nu_U(x,t;y,s)\d x \,=\, 1~, \quad 
  \hbox{and} \quad \inttwo \Gamma^\nu_U (x,t;y,s)\d y \,=\, 1~.
\end{equation}
Note that the first equality uses the fact that $U$ is 
divergence-free. 

\medskip We now consider the particular case where $\omega(x,t) =
\omega^\nu(x,t)$ and $U(x,t) = u^\nu(x,t)$. Then $\partial_t \omega +
(U\cdot\nabla)\omega = \nu \Delta \omega$ by construction, and the
results established in \cite{GMO88} show that assumption \eqref{Ubdd}
is satisfied with $C_0 = C|\alpha|/\nu$, where $C > 0$ is a universal
constant and $|\alpha| \,=\, |\alpha_1| + \dots + |\alpha_N|$.
Using the H\"older continuity \eqref{Gammcont} and the fact that
$\omega^\nu(\cdot,t) \weakto \mu$ as $t \to 0$, we can take the 
limit $s \to 0$ in the representation \eqref{omegrep} and obtain, 
for any $\nu > 0$, the following expression
\begin{equation}\label{omegrep2}
  \omega^\nu(x,t) \,=\, \inttwo \Gamma^\nu_{u^\nu}(x,t;y,0)
  \d \mu_y \,=\, \sum_{i=1}^N \alpha_i\Gamma^\nu_{u^\nu}(x,t;x_i,0)~.
\end{equation}
Setting $\omega_i^\nu(x,t) = \alpha_i\Gamma^\nu_{u^\nu}(x,t;x_i,0)$, 
we obtain the desired decomposition \eqref{omdecomp}, and the 
various properties of $\omega_i^\nu(x,t)$ follow directly from
\eqref{Gammup} and \eqref{Gammint}. \QED

\bigskip\noindent{\bf Proof of Lemma~\ref{PWapprox}.} For simplicity, 
we set $z^\nu(t) = (z^\nu_1(t),\dots,z^\nu_N(t)) \in (\real^2)^N$
and we rewrite system \eqref{PW2} as $\dot z^\nu(t) = 
F(z^\nu(t),\nu t)$, where $F : (\real^2)^N \times (0,\infty) 
\to (\real^2)^N$ is defined by
\begin{equation}\label{FFdef}
  F_i(z,\eta) \,=\, \sum_{j\neq i}\frac{\alpha_j}{\sqrt{\eta}}
  \,v^G\Bigl(\frac{z_i-z_j}{\sqrt{\eta}}\Bigr)~, \qquad 
  i \in \{1,\dots,N\}~.
\end{equation}
For any $\delta \ge 0$, we denote $\Omega_\delta = \{z \in (\real^2)^N
\,;\, |z_i-z_j| > \delta \hbox{ for all }i\neq j\}$. Then $F$ extends
to a smooth map from $\Omega_0 \times [0,\infty)$ to $(\real^2)^N$, 
with
\begin{equation}\label{FF0def}
  F_i(z,0) \,=\, \sum_{j\neq i}\frac{\alpha_j}{2\pi}
  \,\frac{(z_i-z_j)^\perp}{|z_i-z_j|^2}~, \qquad 
  i \in \{1,\dots,N\}~,
\end{equation}
This remark already implies that system \eqref{PW2} is locally
well-posed for all initial data in $\Omega_0$. Global well-posedness
easily follows, because as soon as $\eta$ is bounded away from zero,
the vector field $z \mapsto F(z,\eta)$ is smooth and uniformly
bounded.

We now compare the solutions of \eqref{PW2} in $\Omega_0$ with 
those of system \eqref{PW}, which can be written as $\dot z(t) =
F(z(t),0)$. If $z \in \Omega_\delta$ for some $\delta > 0$, then 
using \eqref{FFdef}, \eqref{FF0def} and the definition \eqref{Gdef} 
of $v^G$, we easily find
\[
  |F_i(z,\eta) - F_i(z,0)| \,\le\, \sum_{j\neq i}\frac{|\alpha_j|}{2\pi}
  \,\frac{1}{|z_i-z_j|}\,e^{-|z_i-z_j|^2/(4\eta)} \,\le\, 
  \frac{|\alpha|}{2\pi \delta}\,e^{-\delta^2/(4\eta)}~,
\]
for any $\eta > 0$. Similarly, if $z, \tilde z \in \Omega_\delta$, then
\[ 
  |F_i(z,0) - F_i(\tilde z,0)| \,\le\, \sum_{j\neq i}\frac{|\alpha_j|}
  {\pi\delta^2}\,\max\{|z_1-\tilde z_1|,\dots,|z_N-\tilde z_N|\} 
  \,\le\, \frac{|\alpha|}{\pi\delta^2}\,\|z-\tilde z\|~,
\]
where $\|z-\tilde z\| = \max\{|z_1-\tilde z_1|,\dots,|z_N-\tilde z_N|\}$. 

Assume now that $z \in C^0([0,T],(\real^2)^N)$ is a solution of 
\eqref{PW} satisfying \eqref{dmin} for some $d > 0$, and take $\delta 
\in (0,d)$. For any $\nu > 0$, let $z^\nu(t)$ denote the unique 
solution of \eqref{PW2} with initial data $z^\nu(0) = z(0)$. As long 
as $z^\nu(t)$ stays in $\Omega_\delta$, we have
\begin{align*}
  \|\dot z^\nu(t) - \dot z(t)\| \,&\le\, \|F(z^\nu(t),\nu t) - 
  F(z^\nu(t),0)\| + \|F(z^\nu(t),0) - F(z(t),0)\| \\
  \,&\le\, \frac{|\alpha|}{2\pi \delta}\,e^{-\delta^2/(4\nu t)} + 
  \frac{|\alpha|}{\pi \delta^2}\,\|z^\nu(t) - z(t)\|~,
\end{align*}
hence
\begin{equation}\label{zcomp}
  \|z^\nu(t) - z(t)\| \,\le\, \frac{|\alpha|}{2\pi \delta}
  \int_0^t e^{|\alpha|(t-s)/(\pi\delta^2)}\,e^{-\delta^2/(4\nu s)}
  \d s \,\le\, \frac{\delta}{2}\, e^{|\alpha|t/(\pi\delta^2)}
  \,e^{-\delta^2/(4\nu t)}~.
\end{equation}
If $\nu > 0$ is sufficiently small, this implies that $z^\nu(t) \in 
\Omega_\delta$ for all $t \in [0,T]$, hence \eqref{zcomp} holds for
$t \in [0,T]$. Choosing for instance $\delta = d\sqrt{4/5}$, we obtain
\eqref{PWcomp} with $K_1 = \exp(CT/T_0)$. For larger values 
of $\nu$, the solution $z^\nu(t)$ may leave $\Omega_\delta$, but in 
that case the bound \eqref{PWcomp} still holds if we take the constant
$K_1$ large enough. \QED

\bigskip\noindent{\bf Proof of Lemma~\ref{Rexpand}.} Let 
$r = |\xi|/|\eta| < 1$, and $\psi = \theta-\phi$. We have
$$
  |\xi + \eta|^2 \,=\, |\eta|^2 (1+2r\cos(\psi)+r^2)
  \,=\, |\eta|^2\,|1+z|^2~,
$$
where $z = r\,e^{i\psi} \in \complex$. Now
\begin{align*} 
  \frac{1}{|1+z|^2} \,&=\, \Bigl(1-z+z^2-z^3+\dots\Bigr)
  \Bigl(1-\bar z + \bar z^2 - \bar z^3+\dots\Bigr) \\
  \,&=\, 1 - (z+\bar z) + (z^2 + z\bar z + \bar z^2) - 
  (z^3 + z^2 \bar z + z\bar z^2 + \bar z^3) + \dots~.
\end{align*} 
But, for each $n \in \natural$, 
\[ 
  z^n + z^{n-1}\bar z + \dots + z\bar z^{n-1} + \bar z^n 
  \,=\, \frac{z^{n+1} - {\bar z}^{n+1}}{z - \bar z}
  \,=\, r^n\,\frac{\sin((n+1)\psi)}{\sin(\psi)}~,
\]
if $\sin(\psi) \neq 0$. Thus
\begin{equation}\label{geom}
  \frac{1}{|\xi+\eta|^2} - \frac{1}{|\eta|^2} \,=\, 
  \frac{1}{|\eta|^2}\Bigl(\frac{1}{|1+z|^2} - 1\Bigr)
  \,=\, \frac{1}{|\eta|^2}\sum_{n=1}^\infty (-1)^n
  \frac{|\xi|^n}{|\eta|^n}\, \,\frac{\sin((n+1)\psi)}{\sin(\psi)}~.
\end{equation}
Multiplying the first and the last member of \eqref{geom} by 
$\xi\cdot\eta^\perp = |\xi||\eta|\sin(\psi)$, we obtain \eqref{V11}.
This concludes the proof. \QED

\medskip\noindent{\bf Acknowledgements.} This work started as
a collaboration with Isabelle Gallagher, to whom I am indebted 
for significant help at the early stage of the project. It is also
a pleasure to thank Maurice Rossi for illuminating discussions 
on the physical signification of the results presented here. 
The author was supported in part by the ANR project ``PREFERED''
of the French Ministry of Research.


\end{document}